\documentclass[12pt,a4paper]{article}

\usepackage[english]{babel}
\usepackage[utf8x]{inputenc}
\usepackage[T1]{fontenc}

\usepackage{geometry}
\usepackage{amsmath,amssymb,amsthm}
\usepackage{graphicx}
\graphicspath{{images/}}
\usepackage[colorinlistoftodos]{todonotes}
\usepackage[colorlinks=true, allcolors=blue]{hyperref}

\usepackage[numbered,framed]{matlab-prettifier}

\usepackage{todonotes}

\usepackage{algorithm}
\usepackage{algpseudocode}

\usepackage{color}

\usepackage{fancyhdr}
\newcommand{\TV}{\mathrm{TV}}
\newcommand{\TVd}{\mathrm{TV}_d}
\newcommand{\C}{\mathbb{C}}
\newcommand{\R}{\mathbb{R}}

\newcommand{\ICBTV}{\mathrm{ICB}_\TV}

\newcommand{\dx}{\,\mathrm{d}x}

\newcommand{\fbold}{\mathbf{f}}
\newcommand{\ubold}{\mathbf{u}}

\newcommand{\ybold}{\mathbf{y}}
\newcommand{\zbold}{\mathbf{z}}
\newcommand{\Tbold}{\mathbf{T}}

\newcommand{\pa}{p_0^\eta}
\newcommand{\qa}{q_0^\eta}

\newcommand{\diverg}{\mathrm{div}}

\newcommand{\Kcal}{\mathcal{K}}

\newcommand{\Fcal}{\mathcal{F}}
\newcommand{\Scal}{\mathcal{S}}
\newcommand{\Pcal}{\mathcal{P}}

\newcommand{\Real}{\mathrm{Re}}
\newcommand{\Imag}{\mathrm{Im}}

\newcommand{\prox}{\mathrm{prox}}
\newcommand{\proj}{\mathrm{proj}}

\numberwithin{equation}{section}

\begin{document}
\title{Dynamic MRI Reconstruction from Undersampled Data with an Anatomical Prescan}
\author{Julian Rasch\thanks{Applied Mathematics M\"unster: Institute for Analysis and Computational Mathematics, 
Westf\"alische Wilhelms-Universit\"at (WWU) M\"unster. Einsteinstr. 62, 48149 M\"unster, Germany. e-mail: julian.rasch@wwu.de}, 
Ville Kolehmainen\thanks{Department of Applied Physics, University of Eastern Finland, POB1627, 70211 Kuopio, Finland},
Riikka Nivaj\"arvi\thanks{Kuopio Biomedical Imaging Unit,
A. I. Virtanen Insititute for Molecular Sciences, 
University of Eastern Finland, POB1627, 70211 Kuopio, Finland}, 
Mikko Kettunen\footnotemark[3], \\
Olli Gr\"ohn\footnotemark[3], 
Martin Burger\footnotemark[1] \ and 
Eva-Maria Brinkmann\footnotemark[1]}

\maketitle
\abstract{
The goal of dynamic magnetic resonance imaging (dynamic MRI) is to visualize tissue properties and their local changes over time that are traceable in the MR signal.
We propose a new variational approach for the reconstruction of subsampled dynamic MR data, which combines smooth, temporal regularization with spatial total variation regularization. In particular, it furthermore uses the infimal convolution of two total variation Bregman distances to incorporate structural a-priori information from an anatomical MRI prescan into the reconstruction of the dynamic image sequence. 
The method promotes the reconstructed image sequence to have a high structural similarity to the anatomical prior, while still allowing for local intensity changes which are smooth in time.
The approach is evaluated using artificial data simulating functional magnetic resonance imaging (fMRI), and experimental dynamic contrast-enhanced magnetic resonance data from small animal imaging using radial golden angle sampling of the $k$-space. \\

\noindent {\bf Keywords: } Dynamic magnetic resonance imaging, spatio-temporal regularization, structural prior, infimal convolution of Bregman distances, total variation, golden angle subsampling
}

\section{Introduction}
\label{sec:intro}
Since its invention in the 1970s, magnetic resonance imaging (MRI) has developed into a very versatile non-invasive medical imaging technique, which can provide information about a rich variety of tissue properties with a high spatial resolution and good soft tissue contrast. 
Initially, however, this high image quality came along with long acquisition times rendering the imaging of rapid physiological changes impossible. 
Hence, since the early days of MRI, the speed-up of the data acquisition process has been a major direction of research. 
In the meantime significant progress in hardware and development of fast imaging protocols has paved the way for modern dynamic imaging techniques such as functional MRI (fMRI) and dynamic contrast-enhanced MRI (DCE-MRI). 

In this paper, we propose a novel variational approach for the reconstruction of highly subsampled dynamic MR data. 
It combines temporal smoothing with spatial total variation (TV) regularization and uses an $\ICBTV$ regularization functional to incorporate structural prior information from an anatomical MR scan acquired prior to the dynamic sequence. 
In the subsequent paragraphs, we briefly introduce the concept of fMRI and DCE-MRI and explain the idea of our approach in more detail afterwards.\\

\noindent {\it Functional magnetic resonance imaging}\\
\noindent
The first fMRI images depicting physiological function inside the brain were presented in the early 1990's \cite{ogawa90,ogawa92}. 
Since then, fMRI has developed into a widely used imaging technique in basic brain research and serves as a tool for the identification and characterization of brain diseases such as Alzheimer's disease.
The fMRI technique relies on the idea that neuronal activity leads to an increase in metabolism and hence to a local change in the cerebral blood flow that can be detected in T$2^*$ weighted MR signals. 
Consequently, during an fMRI experiment a series of fast T$2^*$ weighted MR scans is taken over time to obtain spatial information about the hemodynamic changes in the brain. For a general reference on fMRI, see for example \cite{fmribook}.\\

\noindent {\it Dynamic contrast-enhanced magnetic resonance imaging}\\
\noindent
The technique of DCE-MRI also dates back to the early 1990s \cite{tofts91}. 
In contrast to fMRI, this imaging method is based on injection of an exogenous bolus of gadolinium-based contrast agent.  
Aiming at the detection of cerebral ischemia, brain tumor or infections, it relies on the fact that a characteristic property of normal and healthy tissue in the brain is an intact blood-brain barrier. 
This blood-brain barrier prevents the contrast agent from passing the cell membranes and hence causes the contrast agent to quickly re-enter the vessels instead of accumulating in the tissue. 
However, in tissue that has been damaged by a tumor or an infection the blood-brain barrier is disrupted, and hence the contrast agent accumulates.
In order to detect this change of contrast, a time series of fast MR scans is acquired.
The scans can be based on either a short or long echo time depending on the application. 
While imaging of tumors is usually based on a short echo time, aiming at good T$1$ contrast, the brain ischemia studies usually use T$2^\ast$ weighted imaging with short echo times for detection of the fast perfusion related changes. For further details on DCE-MRI, see for example
the review papers \cite{padhani2002,choyke2003,tofts1999}.
\\ 

\noindent {\it Compressed sensing in magnetic resonance imaging}\\
\noindent
Since the strength of dynamic MRI methods lies in their ability to measure spatial contrast or physiological changes with rapidly switching dynamics, a critical point is to guarantee a sufficiently high temporal resolution in the measured data. 
In this context, the evolution of the theory of compressed sensing (CS) and its application to MRI, e.g. \cite{Candes:Robust,Donoho:CompressedSensing,Lustig:Sparse,Huang:CSinMR}, allowed for a significant acceleration of the data acquisition process. 
Building on the observation that images are likely to have a sparse representation in some transform domain, it led to a boost in the development of new mathematical reconstruction techniques.
These are capable of producing reconstructions of comparable quality from a substantially reduced amount of sampled $k$-space coefficients, provided that the remaining coefficients are chosen in a suitable way. 
In this context, the subsampling scheme is regarded as being suitable, if the resulting aliasing artifacts in the images obtained by linear reconstruction methods are incoherent, or in other words noise-like, in the respective transform domain.
Appropriate nonlinear reconstruction methods promoting sparsity in the particular transform domain and consistency with the acquired data can then facilitate high quality reconstructions despite the comparably small amount of measured data. 
In the following, we will therefore assume that at each timestep of the respective dynamic MRI acquisition process only a very sparse subset of the entire $k$-space is sampled. \\

\noindent {\it The proposed approach}\\
\noindent
Gathering the ideas of fMRI, DCE-MRI and compressed sensing described above, we introduce the three buidling blocks of our variational model for the reconstruction of dynamic MR data. 

The first is an appropriate data fidelity term (cf. Sections \ref{subsec:operator modeling} and \ref{subsec:dynamic imaging}) combined with a (spatial) total variation (TV) regularization \cite{Rudin:ROF} to ensure that the reconstructed images to a certain degree exhibit sparsity in the gradient domain (cf. Section \ref{subsubsec:TV}). 
Provided a sufficiently high temporal resolution, we moreover expect the intensity values to change only in a small part of the entire image domain between two consecutive timesteps (since only a small area is activated). 
Note that this is a common assumption for dynamic MRI applications, especially when motion is negligible, see for example \cite{feng2013,feng2014,wundrak2016}.
Hence, we also suppose a temporal redundancy of the dynamic data that we want to derive benefit from. 
Against this background, it seems reasonable to sample the $k$-space for successive timesteps in a complementary manner, meaning that the subsampling scheme should not be identical among timesteps that are temporally close to each other, but rather such that it adds to the information gained at neighboring points in time. 
One such sampling scheme, which is popular in DCE-MRI, is the radial golden angle sampling where the $k$-space is sampled continuously along radial spokes with a uniform angular stepping ($111.25^\circ$).
It leads to a rather regular spatial coverage of the $k$-space over time with high temporal incoherence and offers great flexibility with respect to the selection of the time resolution in the sense that one can afterwards choose how many radial spokes per frame are used in the reconstruction of the data \cite{winkelmann2007,feng2014}.
The second building block is hence a temporal regularization (cf. Section \ref{subsec:temporal regularization}), which couples the sampled data across the dynamic sequence such that the information acquired at neighboring time points can contribute to the reconstruction of the respective frame.

Finally, we seek to take advantage of the specific experimental setups:
in both fMRI and DCE-MRI, the experimental protocol includes an anatomical high-resolution MRI scan, which is carried out prior to the acquisition of the actual dynamic MRI data. 
While this data contains high resolution anatomical information of the target, it is mostly used for visualization purposes only. 
One of the goals of our approach is to utilize the structural information contained in the anatomical image as prior information in the reconstruction of the dynamic image sequence.
In order to increase quality and hence obtain a higher spatial resolution in the reconstruction of the dynamic MRI data series, we incorporate the edge information of the anatomical high-resolution scan encoded in the subgradient into the reconstruction framework (cf. Section \ref{subsubsec:structural prior}). 
Mathematically, this is realized by the $\ICBTV$ regularization that has first been introduced for color image processing in \cite{Moeller:ColorBregmanTV} and discussed in more detail and applied to PET-MRI joint reconstruction in \cite{Rasch2017}.\\ 

\noindent {\it Related work}\\
\noindent
Our approach certainly stands in the tradition of other CS-related variational models for the reconstruction of dynamic MRI data that combine suitable undersampling strategies with spatio-temporal regularization \cite{Lustig:ktSPARSE,Adluru2007_2,Jung:ktFOCUSS,Tremoulheac:lowRankPlusSparsePrior,Otazo:lowRankPlusSparseMatrixDecomposition,Yao:nuclearNormdMRI,feng2013,feng2014,wundrak2016}. 
Note that the list of related work given here is only a small excerpt of methods proposed in literature. 
We believe that the method of Adluru and coworkers \cite{Adluru2007_2} is closest to our new approach.
However, to the best of our knowledge, no method has yet been proposed which combines the incorporation of high-resolution structural a-priori information from the anatomical prescan with spatio-temporal regularization allowing for both simultaneously, a surprisingly high temporal and spatial resolution.\\

\noindent {\it Organization of the paper}\\
\noindent
The remainder of the paper is organized as follows: In Section 2 we describe the mathematical forward model for MRI and introduce the proposed reconstruction method.  
In Section 3, we first comment on the implementation of all operators and then explain how to solve the minimization problem by a first-order primal-dual optimization method. 
The numerical results are given in Section 4. 
First we evaluate the method using a test case based on fMRI with simulated measurement data, then the approach is evaluated using experimental small animal DCE-MRI data from a glioma model rat specimen, where in both cases the data results from golden angle radial sampling. 
We complete the paper with a brief conclusion and give an outlook on future directions of research in this area. 
\section{The proposed reconstruction framework}
\label{sec:modeling}
In this section we translate our novel reconstruction approach, whose basic modules we motivated above, into a mathematical setting. 
To this end, we specify the modeling of the required operators and all terms contained.
We start with the static MRI forward problem suitable for modeling any arbitrary MR data acquisition protocol. 

\subsection{Undersampled MRI: A mathematical model}
\label{subsec:operator modeling}

While the forward model in MRI can become arbitrarily difficult if modeling the influence of all the underlying physical and technical parameters, it is common practice to assume it to be a simple Fourier transform and hide all the unkowns.
Following this line, we define the Fourier transform of the MR image $u \colon \Omega \to \C$ at $k$-space coordinate $\xi$ by 
\begin{align}\label{eq:fourier_transform_cont}
	(\Kcal u)(\xi) = \int_\Omega u(x) e^{-i x \cdot \xi} \dx, 
\end{align}
where $\Omega \subset \R^d$ is a bounded image domain (usually $d=2,3$).
The associated inverse problem is then to reconstruct $u$ from noisy measurements $f$ obtained by 
\begin{align}\label{eq:inverse_problem}
	\Kcal u + e = f, 
\end{align}
where $e$ models the (complex-valued) noise. 

A strong point of MRI is that it allows for a great variety of different data acquisition procedures each of which leads to a very distinct appearance of the image $u$, which does not become obvious from the above formula \eqref{eq:fourier_transform_cont}, since, as mentioned before, all the parameters characterizing the specific acquisition scheme are hidden.   
In particular, different scanning protocols produce different types of image contrasts.
For example, depending on the respective relaxivities, the same tissue can appear very bright in T1 contrast while being very dark in T2 and some details might only be visible at a specific setup of the acquisition scheme.    
However, since the relaxivity varies only between different types of tissues while being constant across the same tissue, it is very likely for images of different contrast to nevertheless share the same structures (with the exception of features just visible in one of the contrasts).

In practice, one has access to only a small portion of Fourier coefficients $\xi_m$, $m = 1, \dots, M$, such that the system \eqref{eq:inverse_problem} becomes highly underdetermined. 
As a remedy, one also discretizes the image domain $\Omega$ and samples the desired reconstruction $u$ at the respective locations. 
More precisely, letting $\Omega_N$ be a discretization of the image domain $\Omega$ into $N$ pixels, we define the discrete Fourier transform $\Kcal_N \colon \C^N \to \C^M$ of $u$ by 
\begin{align}\label{eq:fourier_transform}
  (\Kcal_N u)(\xi_m) = \sum_{x_n \in \Omega_N} u(x_n) e^{-i x_n \cdot \xi_m} 
\end{align}
for all $k$-space coordinates $\xi_m$, $m=1,\dots,M$.
Here $u(x_n)$ denotes a sampling of $u$ at location $x_n \in \Omega_N$.
It is worth noticing that the $k$-space coordinates $\xi_m$ can be chosen arbitrarily and are in particular independent of the choice of discretization of the image domain.
Hence, from a mathematical perspective, the resolution $N$ of the reconstruction can (theoretically) be chosen arbitrarily at the price of a nontrivial nullspace, and in particular does not depend on the number of sampled Fourier coefficients.
The standard way in a Cartesian setting is to choose $M=N$ , i.e. to sample as many Fourier coefficients as pixels in the desired image, since in this case there holds a one-to-one relationship between an image and its Fourier transform. 
Vice versa, having measured $M$ Fourier coefficients on a Cartesian grid, the canonical resolution of $u$ is $M$ pixels. 
This gives the notion of the standard Cartesian Fourier transform, which can simply be inverted (cf. Section \ref{subsec:implementation of operators}).

However, driven in particular by developments in compressed sensing (cf. \cite{Lustig:Sparse,Lustig:ktSPARSE}) it is common practice to speed up the data acquisition by undersampling the $k$-space, i.e. by performing significantly less measurements $M \ll N$ than desired image pixels. 
In particular, the measurements are usually not sampled on a Cartesian grid (cf. Section \ref{subsec:dynamic imaging}).
In this case there is neither the possibility to simply ``invert'' the operator, nor a ``canonical'' way of choosing the resolution $N$. 
Instead one can try to compensate for the missing data using additional a-priori information (regularization) and choose the resolution as high as the respective method allows. 
The quality of the reconstruction then depends on the location of the measured Fourier coefficients and on the ratio $M/N$. 
As we show in the numerical section, this ratio can actually be less than one per cent (as opposed to $100$ per cent for $M=N$), choosing a golden angle sampling and an appropriate reconstruction method (see also Section \ref{sec:numerics}).

Summing up the discrete setting, we aim to solve the inverse problem
\begin{align}\label{eq:forward_model}
	\Kcal_N u + e = f,
\end{align}
with $u \in \C^N$, $f \in \C^M$ and random measurements errors $e \in \C^M$. 
Since the noise in MRI is commonly approximated by  Gaussian ones (and more advanced noise models such Rician are beyond the scope of the present study), the canonical related minimization problem reads
\begin{align*}
	u \in \arg \min_{u \in \C^N} ~ \frac{\alpha}{2} \| \Kcal_N u - f \|_{\C^M}^2 + J(u).
\end{align*}
The regularizer $J$ is yet to be defined. 
The goal is to construct it such that it promotes desirable properties (such as smoothness etc.) of the sought-after solution.

\subsection{Dynamic MRI: model extension}
\label{subsec:dynamic imaging}

The scanning protocol we consider throughout this paper involves a (densely sampled) prior scan, followed by a (subsampled) dynamic scan. 
In order to cover the dynamic scan we extend the above idea along the same lines to dynamic MR imaging. 

Let us first comment on the characteristics of the dynamic data: Generally, there are many eligible options in how to choose the specific scanning protocol and (sub-) sampling pattern in practice. 
In this work, we consider radial samplings, that is the data is collected on radial spokes through the $k$-space center.
Note that the basic idea does not change for different types of samplings. 
The interesting part about dynamic MRI is that the MR scanner collects the data almost continuously in time, meaning that it measures the radial spokes in $k$-space sequentially one after the other. 
This leaves us free to determine the time resolution of the recorded sequence only during the reconstruction process. 
More precisely, for the reconstruction we divide the entire collection of spokes $\fbold$ into an arbitrary number $T$ of sets of (consecutive) spokes $\fbold = [f_1, \dots, f_T]$. 
For example, if the scanner recorded 100 spokes of data, we can divide the data set into $T = 100$ sets of 1 spoke each, or $T = 10$ sets of 10 spokes each, etc. leading to a very high time resolution (1 spoke/frame) or to a more moderate time resolution (10 spokes/frame), respectively. 
In choosing $T$, we naturally always have to face the trade-off between a high temporal resolution and a sufficient amount of data per frame available for a meaningful reconstruction to build upon. 
In the following, we shall always assume that we have already divided the dynamic data set into $T$ parts. 

In this setting, we then denote the anatomical prior image by $u_0 \in \C^{N_0}$ with the corresponding densely sampled $k$-space data $f_0$, and for $t = 1, \dots, T$ the undersampled data set $\fbold = [f_1, \dots, f_T]$, where each $f_t \in \C^{M_t}$. 
The corresponding sought-after reconstruction is $\ubold = [u_1, \dots, u_T] \in \C^{N \times T}$.
The imaging operators are 
\begin{align*}
      \Kcal_0 \colon \C^{N_0} \to\C^{M_0}, \qquad \Kcal_t \colon \C^{N} \to \C^{M_t},
\end{align*}
for the (densely) sampled anatomical prior and for time index $t$ in the dynamic sequence, respectively. 
Note that we drop the dependence on the resolution $N_0$, $N$ in the following since it will always be clear from the context.

We remark that we keep $N$ fixed for the dynamic scans, implying that the spatial resolution of the $u_t$ stays the same over time. 
Theoretically, however, there also exists the possibility to change the resolution over time.
The corresponding joint reconstruction problem (see also \cite{Ehrhardt2015,Rasch2017}) for the dynamic sequence then reads 
\begin{align*}
      \ubold \in \arg \min_{\ubold \in \C^{N\times T}} ~ \sum_{t=1}^T \frac{\alpha_t}{2} \| \Kcal_t u_t - f_t \|_{\C^{M_t}}^2 + J(\ubold,u_0),
\end{align*}
with a regularizer $J$ that now depends on the different frames $u_t$ and on the anatomical prior image $u_0$ and thus introduces a coupling among all these images. 
Note that this regularizer can, of course, be composed of several terms.

\subsection{Spatial regularization: TV and $\boldsymbol{\ICBTV}$} 
We now specify our choice of the regularizer $J$. 
In the current section, we will concentrate on penalty terms that do not take into account that we actually aim at reconstructing an entire sequence of consecutive time frames $\ubold$. 
Instead, these terms rather act on each frame $u_t$ independently. 
The question of how we can link the time frames during the reconstruction will be addressed in the subsequent subsection.

\label{subsec:spatial regularization}
\subsubsection{TV: enforcing sparsity in the gradient domain}
\label{subsubsec:TV} 
As motivated and explained before, we assume that at each time step of the dynamic sequence only a small portion of all Fourier coefficients is sampled along radial spokes through the $k$-space center. 
According to the theory of compressed sensing \cite{Candes:Robust,Donoho:CompressedSensing,Lustig:Sparse,Huang:CSinMR}, this requires the data fidelity term to be accompanied by a term that guarantees sparsity in some transform domain.
One popular choice for such a term and the approach we pursue in this paper is the (spatial) total variation regularization.  
First introduced for image denoising in 1992 \cite{Rudin:ROF}, the distinctive feature of this regularization is the promotion of piecewise constant solutions with sharp edges and thus of sparsity with respect to the image gradient. 
For real-valued images $z \in \mathbb{R}^N$, the discrete, isotropic form of TV reads:
\begin{align*}
	\TVd (z) := \| \nabla z \|_1 :=\sum_{n=1}^N |(\nabla z)_n | = \sum_{n=1}^N \sqrt{|(\nabla z)_{n,1}|^2 + |(\nabla z)_{n,2}|^2}.
\end{align*}
with $| \cdot |$ denoting the Euclidean norm on $\mathbb{R}^N$.

However, since in this work we model the unknown MRI reconstruction to be complex-valued, it is necessary to extend the concept of gradients and total variation accordingly. 
The extension of gradients is straightforward and for a complex-valued image $u \in \C^N$ we just define the gradient by $\nabla \colon \C^N \to \C^{N \times 2}$, where more details on our particular implementation can be found in Section \ref{subsec:implementation of operators}. 
The extension of the total variation is slightly more involved: 
though the underlying image $u$ is complex-valued, in most real-world applications one is interested only in the magnitude of the reconstructed image. 
More precisely, using the identity $u = r e^{i \varphi}$, one ideally aims to perform the regularization of the sought-after image directly on its magnitude, i.e. to penalize the total variation of $r$. 
Unfortunately, this leads to a nonlinear (in $\varphi$) forward operator $\tilde{K} (r,\varphi)$ in \eqref{eq:inverse_problem} or \eqref{eq:forward_model}, rendering the numerical solution of the resulting minimization problem considerably more difficult.   
We refer the reader to \cite{Valkonen2014} for an extensive study of this situation.

In this work, we shall instead use the linear forward model \eqref{eq:forward_model} and an approximation to the total variation on the magnitude, which turns out to have a similar effect in practice. 
More precisely, in the complex-valued case, we redefine the discrete, isotropic total variation $\TVd \colon \C^N \to [0,\infty)$ in the following way 
\begin{align*}
	\TVd (u) = \sum_{n=1}^N \sqrt{|\Real (\nabla u)_{n,1}|^2 + |\Imag (\nabla u)_{n,1}|^2 + |\Real (\nabla u)_{n,2}|^2 + |\Imag (\nabla u)_{n,2}|^2},
\end{align*}
which reveals that this approach is equivalent to regarding the image $u$ to be located in $\R^{N \times 2}$ and using real-valued isotropic total variation on its gradient in $\R^{N \times 4}$.
Simply speaking, this approach consists in taking the magnitude of the gradient of $u$ rather than the gradient of the magnitude of $u$.

At the end of this paragraph we address the TV regularization also in a slightly different context, namely with respect to the reconstruction of the anatomical prior $u_0$. 
As already explained in Section \ref{sec:intro}, we shall use the structural information of this high resolution image as additional a-priori information for the reconstruction of the dynamic sequence $\ubold$.
In order to reconstruct the usually fully sampled prior scan, i.e. M = N in Equation \eqref{eq:fourier_transform}, one would typically simply apply the inverse Fourier transform. 
However, since we do not want to incorporate minor oscillations, but rather like to concentrate on major structures encoded in the prior scan, we incorporate a TV regularization term into the reconstruction problem, i.e., in the first step we solve:
\begin{equation*} \label{anatomicTV}
u_0 \in \arg \min_{u_0 \in \C^{N_0}} ~ \frac{\alpha_0}{2} \| \Kcal_0 u_0 - f_0 \|_{\C^{M_0}}^2 + \TVd(u_0).
\end{equation*}
We hence obtain a high-quality piecewise constant anatomical reconstruction $u_0$ with sharp edges which serves as the basis for our further considerations.

\subsubsection{$\boldsymbol{\ICBTV}$: incorporating structural prior information}
\label{subsubsec:structural prior}
In this paragraph, we explain how we can utilize the structural information of the (high resolution) MR image $u_0$ to support the reconstruction of the dynamic MRI sequence.
For this purpose, we will again build upon the concept of total variation regularization \cite{Rudin:ROF} and the related Bregman distances and iterations \cite{Osher:AnIterativeRegularizationMethod}. 
Extensions of these techniques have been shown to be very effective for coupling the edge information of different color channels in RGB image processing \cite{Moeller:ColorBregmanTV} and of different medical imaging modalities in PET-MRI \cite{Rasch2017}, as well as in the context of related debiasing methods \cite{2016debiasing}.

A natural approach for extracting the structural information from the reconstruction of the prior scan is to use the image gradient, since it essentially encodes the edges and hence the structure of the image. Based on this idea, there emerged a variety of methods for structural priors in the literature, of which we only mention a few here \cite{kaipio1999,Ehrhardt2015,Ehrhardt2016,Bowsher:Bayesian,Chan:Regularized, Nuyts:Mutual,Atre:Evaluation, Leahy1991, Lipinski:Expectation, Vunckx:Evaluation}.  
However, since most of these approaches are not normalized with respect to the size of the gradient, a common difficulty is to deal with images in different scales, i.e. images sharing the same structure, but having very different intensity ranges and hence different sizes of jumps. 
Considering MRI this problem is of particular importance, since the absolute intensities of the reconstructions can be very different depending on the acquisition protocol chosen.

In \cite{Rasch2017} it has therefore been argued for using a subgradient rather than a gradient to include edge information without considering their size. 
In order to define a straightforward extension of subdifferentials to the complex setting, we need to equip $\C^N$ with the following inner product:
\begin{align*}
\langle u,v \rangle_{\C^N} := \Real (u^* v) \qquad u,v \in \C^N 
\end{align*}
with $^*$ denoting the convex conjugation. 
This inner product induces the same norm as the standard inner product on $\C^N$, but is real-valued, which enables us to consider angles between complex vectors. 
Indeed, for $u,v \in \C^N$ we have 
\begin{align*}
	\cos(\varphi) = \frac{\langle u, v \rangle_{\C^N}}{\|u\|_{\C^N} \|v\|_{\C^N}}, 
\end{align*}
where $\varphi$ denotes the angle between $u$ and $v$. 

We now define the subdifferential of a convex functional $J \colon \C^N \to [0,\infty)$ by 
\begin{align*}
      \partial J(u_t) = \{ p \in \C^N ~|~  J(v) \geq J(u_t) + \langle p, v - u_t \rangle_{\C^N} \; \forall v \in \C^N \}. 
\end{align*}
By the chain rule, a subgradient $p_0$ of the total variation at $u_0$ is then given by 
\begin{align*}
	p_0 \in \partial \TVd (u_0) \Leftrightarrow p_0 = - \diverg (q_0), \quad q_0 \in \partial \|\nabla u_0 \|_1,
\end{align*}
with $q_0 \in \C^{N \times 2}$ such that 
\begin{align*}
	|q_0| \leq 1 \quad \text{ and } \quad q_0 = \frac{\nabla u_0}{|\nabla u_0|} \text{ if } \nabla u_0 \neq 0.
\end{align*}
Hence, we see that where the gradient of $u_0$ does not vanish, the vector field $q_0$ equals the direction of the edge in $u_0$ at this location and this information is encoded in the subgradient. 
At positions with zero gradient, the prior $u_0$ does not provide any structural information (apart from being constant), and $q_0$ is an arbitrary vector from the unit ball. 
It has been argued in \cite{Osher:AnIterativeRegularizationMethod,Moeller:ColorBregmanTV,Rasch2017}, that specific choices of subgradients can provide ''hints`` about structure also at these locations. 
In this work, however, we want to use only the location of edges already clearly visible in the reconstructed anatomical prior.  
Accordingly, we will use the following vector field 
\begin{align*}
	q_0 = \begin{cases}
			\frac{\nabla u_0}{|\nabla u_0|}, &\text{ if } \nabla u_0 \neq 0, \\
            0, &\text{ else,}
		\end{cases}
\end{align*}
to obtain a subgradient $p_0 = -\diverg (q_0) \in \partial \TVd (u_0)$.
The associated Bregman distance can be written as 
\begin{align*}
	D_{\TVd}^{p_0}(u_t,u_0) 
    &= \TVd (u_t) - \langle p_0,u_t \rangle_{\C^N}
    = \| \nabla u_t \|_1 - \langle q_0, \nabla u_t \rangle_{\C^{N \times 2}} \\ 
    &= \sum_{ \{\nabla u_0 \neq 0\} } |(\nabla u_t)_n | \left( 1 - \frac{\langle (\nabla u_t)_n, (\nabla u_0)_n \rangle}{|(\nabla u_t)_n||(\nabla u_0)_n|} \right) + \sum_{\{\nabla u_0 = 0\}} |(\nabla u_t)_n| \\
    &= \sum_{\{\nabla u_0 \neq 0 \} } |(\nabla u_t)_n| (1 - \cos (\varphi_n)  ) + \sum_{\{\nabla u_0 = 0 \} } |(\nabla u_t)_n|,
\end{align*}
where $\varphi_n$ denotes the angle between the vectors. 
Hence we conclude that on the support of $\nabla u_0$ the Bregman distance penalizes the gradient of $u_t$ weighted by its directional deviation from $\frac{\nabla u_0}{|\nabla u_0 |}$, while outside the support of $\nabla u_0$ we obtain a standard TV regularization. 
Consequently, this functional favors aligned edges between $u_t$ and the prior $u_0$ without penalizing the size of the gradients. 
More precisely, the height of the jumps of shared edges can be determined only in dependence on the data term, since once an edge of $u$ is aligned to the edge of $u_0$, i.e., $\cos(\varphi_n) = 1$, the regularization functional vanishes at this position. 
We refer to \cite{Moeller:ColorBregmanTV,Rasch2017} for an extensive discussion of the behavior.

In view of noisy data and with respect to practical applications it often makes sense to use a regularized approximation of the vector field $q_0$: 
\begin{align}\label{eq:subgrad_eta}
	\qa = \begin{cases}
			\frac{\nabla u_0}{|\nabla u_0|}, &\text{ if } |\nabla u_0| \geq \eta, \\
            0, &\text{ else.}
		\end{cases}
\end{align}
In this formulation, the threshold parameter $\eta$ defines a minimum height such that we consider a jump in $u_0$ to be an edge and in this way avoids to falsely identify small oscillations as edges (see also \cite{Ehrhardt2016} for a similar approach). 
In this case, the Bregman distance reads 
\begin{align}
	D_{\TVd}^{\pa}(u_t,u_0) = \sum_{\{|\nabla u_0| \geq \eta \} } |(\nabla u_t)_n| (1 - \cos (\varphi_n)  ) + \sum_{\{|\nabla u_0| < \eta \} } |(\nabla u_t)_n|. 
    \label{eq:Bregman Distance with respect to p0eta}
\end{align}

Due to potentially different image contrasts in the prior scan and the dynamic sequence, it is unfortunately unrewarding to directly use this Bregman distance to include the structural prior information in the reconstruction framework for the dynamic sequence.
To illustrate that this is indeed true, we consider a situation, where taking a step up across the edges of the structure in the prior corresponds to taking a step down across the same structure in the dynamic sequence. 
In this situation, the angle $\varphi_n$ between the two vectors in Equation \eqref{eq:Bregman Distance with respect to p0eta} is close to $\pi$ such that the Bregman distance \eqref{eq:Bregman Distance with respect to p0eta} would highly penalize the edge in $u_t$ at the respective position, which in our setting seems counterintuitive. 
We therefore face a situation, where we expect edge positions to coincide, but where we not necessarily aim at aligning edges to the same direction, but rather would like to adapt the regularizer such that it promotes also the opposing direction of the jump. 
In order to obtain a symmetric measure with respect to vector orientation, we again follow \cite{Moeller:ColorBregmanTV} and \cite{Rasch2017} and employ the infimal convolution of two Bregman distances with respect to subgradients with opposing signs 
\begin{align*}
	\ICBTV^{\pa}(u_t,u_0) := \inf_{\phi + \psi = u_t} ~ D_{\TVd}^{\pa}(\phi,u_0) + D_{\TVd}^{-\pa}(\psi,-u_0).
\end{align*}
This functional operation yields a decomposition of the image $u_t$ into two parts of which one matches the edge set of $u_0$, and the other one matches the edge set of $-u_0$, hence with ''inverted`` contrast. 
Acting as an edge indicator independent of the sign of jumps, the $\ICBTV$ functional thus allows to incorporate the structure of the prior image $u_0$ as a-priori information into the reconstruction of the dynamic sequence. For a more detailed discussion of the behavior of the functional we again refer to \cite{Rasch2017}.

\subsection{Temporal regularization}\label{subsec:temporal regularization}
So far, we have only discussed spatial regularization techniques, which affect all the time frames independently from each other, thus not imposing any particular relation between frames at successive points in time. 
The last ingredient we want to add to our dynamic reconstruction framework is an additional temporal regularization, which interlinks consecutive time frames.

In the literature a whole zoo of approaches can be found, reaching from Tikhonov-type regularization of the time derivative \cite{Adluru2007,Adluru2007_2} via sparsity in some known basis, i.e., \cite{Lustig:ktSPARSE} after applying e.g. the temporal Fourier transform \cite{Tremoulheac:lowRankPlusSparsePrior}, the time derivative \cite{feng2014,wundrak2016} or the discrete cosine transform \cite{fang2016high}, to low rank regularization \cite{Yao:nuclearNormdMRI}. 
In addition, it has been proposed to combine the last two classes to get a decomposition into a low rank part and into a part which is sparse in a known transform domain \cite{Tremoulheac:lowRankPlusSparsePrior,Otazo:lowRankPlusSparseMatrixDecomposition}.
To make full use of the temporal redundancy present in the dynamic data set, one should always choose the regularization in dependence on the expected dynamic behavior. 
For example, for a recurrent dynamic such as the beat of a heart a low rank regularization stands to reason, while in case of objects quickly moving from one region of the image to another sparsity of the time derivative seems to be a more natural choice.
Since the (pixelwise) temporal behavior we expect in our applications is commonly modeled by a smooth curve, we penalize the squared time derivative. 

For a stack of images $\ubold = [u_1, \dots, u_T] \in \C^{N \times T}$, a discrete time derivative can be defined by forward differences, i.e.
\begin{align*}
	\partial_t u_t = \begin{cases}
    	(u_{t+1} - u_t)/(\Delta t), & t = 1, \dots, T-1, \\
        0, & t = T,
	\end{cases} 
\end{align*}
with step size $\Delta t$.
Using a Tikhonov-type regularization for the time derivative we arrive at the penalty function 
\begin{align*}
	\Delta t \sum_{t = 1}^{T-1} \| \partial_t u_t \|_{\C^N}^2 = \frac{1}{\Delta t} \sum_{t = 1}^{T-1}  \| u_{t+1} - u_t \|_{\C^N}^2.
\end{align*}

\subsection{The proposed method}
We are finally able to combine the above considerations to a full method for dynamic MRI with a structural prior. 
The goal of the method is that for all $t = 1, \dots, T$ the reconstructed time frame $u_t$ 
\begin{enumerate}
    \item matches the data $f_t$, i.e. $\| \Kcal_t u_t - f_t \|^2_{\C^{M_t}}$ is small, 
    \item 
     has sparse (spatial) gradient
    , i.e. $\TV(u_t)$ is small, 
    \item has a structure similar to the prior $u_0$, i.e. $\ICBTV^{p_0}(u_t,u_0)$ is small, 
    \item should not be too different from the previous and consecutive frame, i.e. the time derivative $\| \partial_t u_t \|_{\C^N}^2$ is small.
\end{enumerate}
In order to weight between these four a-priori assumptions on $u_t$, we introduce weighting factors $\alpha_t >0, \gamma_t >0 , w_t \in [0,1]$ and the reconstruction method becomes:
 \begin{alignat}{4}
   &\ubold \in \arg &&\min_{\ubold} && \sum_{t=1}^T \frac{\alpha_t}{2} \| \Kcal_t u_t - f_t\|_{\C^{M_t}}^2 \qquad &&\text{(data fidelity)} \notag\\ 
   & && +&& \sum_{t=1}^T w_t \TV_d(u_t) \qquad &&\text{(TV regularization)} \notag \\
   & && +&& \sum_{t=1}^T (1-w_t) \ICBTV^{p_0}(u_t, u_0 ) \qquad &&\text{(structural prior)} \notag \\
   & && +&& \sum_{t=1}^{T-1} \frac{\gamma_t}{2} \| u_{t+1} - u_t \|_{\C^N}^2. \qquad &&\text{(temporal smoothing)} 
\label{eq:dynamic_recon}
\end{alignat}
We quickly outline the behavior. 
The parameters $\alpha_t$ control the data fidelity, introducing data accuracy for the reconstruction. 
The parameters $w_t \in [0,1]$ weight between total variation regularization and proximity to the structural prior. 
The weights are chosen between 0 and 1 since both the $\TV$ and the $\ICBTV$ term enforce spatial regularity. 
For small $w_t$, the reconstructions will have a structure very similar to the prior, for $w_t$ close to 1 the prior does not substantially influence the reconstructions. 
Since we expect dynamics only in some part of the image while the rest stays constant, it is reasonable to keep more (but not all the) weight on the prior. 
The parameters $\gamma_t$ control the temporal smoothness of the reconstructed sequence. 
For the sake of clarity we embed the temporal resolution $\Delta t$ in $\gamma_t$. 
However, switching to a different time resolution, the parameters $\gamma_t$ have to be adjusted related to the change in $\Delta t$.

\section{Numerical implementation and solution}
In this section, we explain how to implement and solve the minimization problem \eqref{eq:dynamic_recon} numerically which, depending on the amount of time steps $T$, can be challenging. 
We derive a general primal-dual algorithm for its solution, before we line out some strategies to reduce computational costs and speed up the implementation at the end of this section. 

\subsection{Gradients and sampling operators}
\label{subsec:implementation of operators}

In order to use the discrete total variation already defined in Section \ref{subsubsec:TV}, we need a discrete gradient operator that maps an image $u \in \C^N$ to its gradient $\nabla u \in \C^{N \times 2}$. 
Following \cite{ChambollePock}, we implement the gradient by standard forward differences. Moreover, we  will use its discrete adjoint, the negative divergence $-\diverg$, defined by the identity $\langle \nabla u, w \rangle_{\C^{N \times 2}} = - \langle u, \diverg(w) \rangle_{\C^N}$.
The inner product on the gradient space $\C^{N\times 2}$ is defined in a straightforward way as 
\begin{align*}
	\langle v,w \rangle_{\C^{N \times 2}} = \Real(v_1^* w_1) + \Real(v_2^* w_2),
\end{align*}
for $v,w \in \C^{N \times 2}$.
For $v \in \C^{N \times 2}$, the (isotropic) 1-norm is defined by 
\begin{align*}
	\|v \|_1 := \sum_{i=1}^N \sqrt{|(v_i)_1|^2 + |(v_i)_2|^2},
\end{align*}
and accordingly the dual $\infty$-norm for $w \in \C^{N \times 2}$ is given by 
\begin{align*}
	\| w \|_\infty := \max_{i=1,\cdots,N} |w_i| = \max_{i=1,\cdots,N} ~ \sqrt{|(w_i)_1|^2 + |(w_i)_2|^2}.
\end{align*}

The sampling operators $\Kcal_t \colon \C^N \to \C^{M_t}$ (and analogously for $\Kcal_0 \colon \C^{N_0} \to \C^{M_0}$) we consider are either a standard fast Fourier transform (FFT) on a Cartesian grid, followed by a projection onto the sampled frequencies, or a non-uniform fast Fourier transform (NUFFT), in case the sampled frequencies are not located on a Cartesian grid \cite{Fessler:NUFFT}. \\

\noindent {\it Fourier transform on a Cartesian grid - the simulated data case}\\
\noindent
In the numerical study on artificial data we use a simple version of the Fourier transform and sampling operator (the same can also be found in e.g. \cite{Ehrhardt2016,Rasch2017}). 
We discretize the image domain on the unit square using an (equi-spaced) Cartesian grid with $N_1 \times N_2$ pixels such that the discrete grid points are given by 
\begin{align*}
 \Omega_N = \left\{ \left(\frac{n_1}{N_1-1}, \frac{n_2}{N_2-1} \right) ~\Big|~ n_1 = 0, \dots, N_1-1, \; n_2 = 0, \dots N_2-1 \right\}.
\end{align*}
We proceed analogously with the $k$-space, i.e. the location of the $(m_1,m_2)$-th Fourier coefficient is given by $(m_1/(N_1-1), m_2/(N_2-1))$. Then, we arrive at the following formula for the (standard) Fourier transform $\Fcal$ applied to $u \in \C^{N_1 \times N_2}$:
\begin{align*}
 (\Fcal u)_{m_1,m_2} = \frac{1}{N_1 N_2} \sum_{n_1=0}^{N_1-1} \sum_{n_2=0}^{N_2-1} u_{n_1,n_2} e^{-2\pi i \left(\frac{n_1 m_1}{N_1} + \frac{n_2 m_2}{N_2} \right)},
\end{align*}
where $ m_1 = 0, \dots, N_1-1, m_2 = 0, \dots, N_2 -1$.
For simplicity, we use a vectorized version such that $\Fcal \colon \C^N \to \C^N$ with $N = N_1 \cdot N_2$.
We then employ a simple sampling operator $\Scal_t \colon \C^N \to \C^{M_t}$ which discards all Fourier frequencies which are not located on the desired sampling geometry at time $t$ (i.e. the chosen spokes). 
More precisely, following \cite{Ehrhardt2016}, if we let $\Pcal_t \colon \{1,\dots,M_t \} \to \{1,\dots,N\}$ be an injective mapping which chooses $M_t$ Fourier coefficients from the $N$ coefficients available, we can define the sampling operator $\Scal$ applied to $f \in \C^N$ as
\begin{align*}
 \Scal_t \colon \C^N \to \C^{M_t}, \quad (\Scal_t f)_k = f_{\Pcal_t(k)}.
\end{align*}
The full forward operator $\Kcal_t$ can hence be expressed as 
\begin{align}\label{eq:forward_op_art}
 \Kcal_t \colon \C^N \xrightarrow{\Fcal} \C^N \xrightarrow{\Scal_t} \C^{M_t}.
\end{align}
The corresponding adjoint operator of $\Kcal_t$ is given by  
\begin{align*}
 \Kcal_t^* \colon \C^{M_t} \xrightarrow{\Scal_t^*} \C^N \xrightarrow{\Fcal^{-1}} \C^N,
\end{align*}
where $\Fcal^{-1}$ denotes the standard inverse Fourier transform and $\Scal_t^*$ `fills' the missing frequencies with zeros, i.e. 
\begin{align*}
 (\Scal_t^*z)_l = \sum_{k=1}^{M_t} z_k \delta_{l,\Pcal_t(k)}, \qquad \text{for } l = 1, \dots, N.
\end{align*}
For the prior $u_0$, we choose a full Cartesian sampling, which corresponds to $\Pcal_0$ being the identity. For the subsequent dynamic scan, we set up $\Pcal_t$ such that it chooses the frequencies located on (discrete) spokes through the center of the $k$-space.     
It is important to notice that this implies that the locations of the (discretized) spokes are still located on a Cartesian grid, which allows to employ a standard fast Fourier transform (FFT) followed by the above projection onto the desired frequencies. 
This is not the case for the operators we use for real data. \\

\noindent {\it Non-uniform Fourier transform - the real data case}\\ 
\noindent 
In contrast to the above (simplified) setup for artificial data, in many real world application the measured $k$-space frequencies $\xi_m$ in \eqref{eq:fourier_transform} are {\it not} located on a Cartesian grid. 
While this is not a problem with respect to the formula itself, it however excludes the possibility to employ a fast Fourier transform, 
which usually reduces the computational costs of an $N$-point Fourier transform from an order of $O(N^2)$ to $O(N \log N)$.
To get to a similar order of convergence also for non-Cartesian samplings, it is necessary to employ the concept of non-uniform fast Fourier transforms (NUFFT) \cite{Fessler:NUFFT,Fessler:code,Matej2004,Nguyen:1999,Strohmer2000}. 
We only give a quick intuition here and for further information we refer the reader to the literature listed above. 
The main idea is to use a (weighted) and oversampled standard Cartesian $K$-point FFT $\Fcal$, $K \geq N$ followed by an interpolation $\Scal$ in $k$-space onto the desired frequencies $\xi_m$. 
Note that the oversampling takes place in $k$-space.
The operator $\Kcal_t$ for time $t$ can hence again be expressed as a concatenation of a $K$-point FFT and a sampling operator 
\textbf{\begin{align*}
 \Kcal_t \colon \C^N \xrightarrow{\Fcal} \C^N \xrightarrow{\Scal_t} \C^{M_t}.
\end{align*}}
For our numerical experiments with the experimental DCE-MRI data, the sampling operator $\Scal_t$ and its adjoint were taken from the NUFFT package \cite{Fessler:code}.

\subsection{Numerical solution}
Due to the nondifferentiablity and the involved operators we apply a primal-dual method \cite{ChambollePock} to solve the minimization problem \eqref{eq:dynamic_recon}. 
We first line out how to solve the (simple) TV-regularized problem for the prior (\ref{tvu0}) and then extend the approach to the dynamic problem. 
Interestingly, the problem for the prior already provides all the ingredients needed for the numerical solution of the dynamic problem, which can then be done in a very straightforward way. 
We consider the problem 
\begin{equation} \label{tvu0}
	\min_{u_0} ~ \frac{\alpha_0}{2} \| \Kcal_0 u_0 - f_0 \|_{\C^{M_0}}^2 + \| \nabla u_0 \|_1, 
\end{equation}
with $u_0 \in \C^{N_0}$.
Dualizing both terms leads to its primal-dual formulation 
\begin{align}\label{eq:tv_pd}
	\min_{u_0} \max_{y_1,y_2} ~ \langle y_1, \Kcal_0 u_0 - f_0 \rangle_{C^{M_0}} - \frac{1}{2 \alpha_0} \|y_1 \|_{\C^{M_0}}^2 + \langle y_2, \nabla u_0 \rangle_{\C^{N_0 \times 2}} + \chi_{C}(y_2),
\end{align}
where $y_1 \in \C^{M_0}$ and $\chi_{C}$ denotes the characteristic function of the set  
\begin{align*}
	C := \{ y \in \C^{N_0 \times 2} ~|~ \|y \|_\infty \leq 1 \}.
\end{align*}
The primal-dual algorithm in \cite{ChambollePock} now essentially consists in performing a proximal gradient descent on the primal variable $u_0$ and a proximal gradient ascent on the dual variables $y_1$ and $y_2$, where the gradients are taken with respect to the linear part, the proximum with respect to the nonlinear part. 
We hence need to compute the proximal operators for the nonlinear parts in \eqref{eq:tv_pd} to obtain the update steps for $u_0$ and $y_1,y_2$. 
It is easy to see that the proximal operator for $\phi(y_1) = \frac{1}{2 \alpha} \|y_1\|_{\C^{M_0}}^2$ is given by 
\begin{align}\label{eq:prox_dual_l2}
	y_1 = \prox_{\sigma \phi} (r) \Leftrightarrow y_1 = \frac{\alpha r}{\alpha + \sigma}.
\end{align}
The proximal operators for the update of $y_2$ are given by a simple projection onto the set $C$, i.e. 
\begin{align}\label{eq:prox_proj}
	y_2 = \proj_C (r) \Leftrightarrow (y_2)_i = r_i / \max (|r_i|,1) \quad \text{for all } i.
\end{align}
Putting everything together leads to Algorithm \ref{alg:prior}.
\begin{algorithm}[t!] 
\caption{\textbf{Reconstruction of the prior}}
{
\begin{algorithmic}[1]
\Require step sizes $\tau,\sigma > 0$, data $f_0$, parameter $\alpha_0$
\Ensure $u_0^0 = \bar{u}_0^0 = \Kcal_0^*f_0, ~ y_1^0 = y_2^0 = 0$
	\While{$\sim$ stop crit}
    	\State {\it Dual updates}
        \State $y_1^{k+1} = (\alpha_0 \left[ y_1^k + \sigma (\Kcal_0 \bar{u}_0^k - f_0)\right]) / (\alpha_0 + \sigma)$
          \State $y_2^{k+1} = \proj_C \left(y_2^k + \sigma \nabla \bar{u}_0^k\right)$
          \State {\it Primal updates}
          \State $u_0^{k+1} =  u_0^k - \tau \left[ \Kcal_0^* y_1^{k+1} - \diverg(y_2^{k+1}) \right]$
          \State {\it Overrelaxation}
          \State $\bar{u}_0^{k+1}= 2 u_0^{k+1} - u_0^k$
	\EndWhile\\
\Return $u_0 = u_0^k$
\end{algorithmic}
}
\label{alg:prior}
\end{algorithm}
\ \\

The numerical realization of the dynamic problem is now straightforward.
In order to deal with the infimal convolution, we use its definition and introduce an additional auxiliary variable yielding 
\begin{alignat*}{4}
	&\min_{\ubold}&& ~ &&\sum_{t=1}^T \frac{\alpha_t}{2} \| \Kcal_t u_t - f_t \|_{\C^{M_t}}^2 + \sum_{t=1}^{T-1} \frac{\gamma_t}{2} \|u_{t+1} - u_t \|_{\C_N}^2 + \sum_{t=1}^T w_t \TV(u_t) \\
	& && + && \sum_{t=1}^T (1-w_t) \ICBTV^{p_0}(u_t,u_0) \\    
   = &\min_{\ubold,\zbold}&& ~ &&\sum_{t=1}^T \frac{\alpha_t}{2} \| \Kcal_t u_t - f_t \|_{\C^{M_t}}^2 + \sum_{t=1}^{T-1} \frac{\gamma_t}{2} \|u_{t+1} - u_t \|_{\C_N}^2 +\sum_{t=1}^T w_t \| \nabla u_t \|_1\\
   & && + &&\sum_{t=1}^T (1-w_t) \left[ \| \nabla (u_t - z_t) \|_1 + \| \nabla z_t \|_1  - \langle p_0,u_t \rangle_{\C^N} + \langle 2 p_0, z_t \rangle_{\C^N} \right]
\end{alignat*}
where $\ubold = [u_1, \dots, u_T] \in \C^{N \times T}$ and $\zbold = [z_1, \dots, z_T] \in \C^{N \times T}$.
Introducing a dual variable $\ybold$ for all the terms containing an operator, leads to the primal-dual formulation 
\begin{alignat}{4}
\label{eq:primal_dual}
	&\min_{\ubold,\zbold} \max_{\ybold} && ~ && \sum_{t=1}^T \left(\langle y_{t,1}, \Kcal_t u_t - f_t \rangle_{\C^{M_t}} - \frac{1}{2 \alpha_t} \|y_{t,1} \|_{\C^M}^2 \right) + \sum_{t=1}^{T-1} \frac{\gamma_t}{2} \| u_{t+1} - u_t \|_{\C^N}^2 \notag \\
    & && + && \sum_{t=1}^T \left(\langle y_{t,2}, \nabla u_t \rangle_{\C^{N \times 2}} + \langle y_{t,3}, \nabla (u_t - z_t) _{\C^{N \times 2}} + \langle y_{t,4}, \nabla z_t \rangle_{\C^{N \times 2}}\right) \notag \\
    & && - && \sum_{t=1}^T \langle (1-w_t)p_0,u_t \rangle_{\C^N} + \sum_{t=1}^T \langle 2(1-w_t)p_0,z_t \rangle_{\C^N}  \notag \\
    & && + && \sum_{t=1}^T \left(\chi_{C_{t,2}}(y_{t,2}) + \chi_{C_{t,3}}(y_{t,3}) + \chi_{C_{t,4}}(y_{t,4})\right)
\end{alignat}
where $\ybold = [\ybold_1, \dots, \ybold_T]$, $\ybold_t = [y_{t,1}, \dots, y_{t,4}]$, and for all $t = 1, \dots, T$, $u_{t} \in \C^{M_t}$ and
\begin{align*}
	&C_{t,2} := \{ y \in \C^{M \times 2} ~|~ \|y \|_{\infty} \leq w_t \}, \\
    &C_{t,3} := \{ y \in \C^{M \times 2} ~|~ \|y \|_{\infty} \leq (1-w_t) \}, \\
    &C_{t,4} := \{ y \in \C^{M \times 2} ~|~ \|y \|_{\infty} \leq (1-w_t) \}. \\
\end{align*}
\begin{algorithm}[t!]
\caption{\textbf{Dynamic reconstruction with structural prior}}
{
\begin{algorithmic}[1]
\Require step sizes $\tau,\sigma > 0$, subgradient $p_0$, for all $t=1,\dots,T$: data $f_t$, parameters $\alpha_t$, $w_t$, $\gamma_t$
\Ensure for all $t=1,\dots,T$: $u_t^0 = \bar{u}_t^0 = \Kcal_t^*f_t, ~ z_t^0 = \bar{z}_t^0 = 0, ~ y_{t,1}^0 = y_{t,2}^0 = y_{t,3}^0 = y_{t,4}^0 = 0$
	\While{$\sim$ stop crit}
    	\For{t=1,\dots,T} 
          \State {\it Dual updates}
          \State $y_{t,1}^{k+1} = \frac{\alpha_t \left[ y_{t,1} + \sigma (\Kcal_t \bar{u}_t^k - f_t)\right]}{\alpha_t + \sigma}$
          \State $y_{t,2}^{k+1} = \proj_{C_2}\left(y_{t,2}^k + \sigma \nabla \bar{u}_t^k\right)$
          \State $y_{t,3}^{k+1} = \proj_{C_3}\left(y_{t,3}^k + \sigma \nabla (\bar{u}_t^k - \bar{z}_t^k) \right)$
          \State $y_{t,4}^{k+1} = \proj_{C_4}\left(y_{t,4}^k + \sigma \nabla \bar{z}_t^k \right)$
          \State {\it Primal updates}
          \State $u_t^{k+1} =  \frac{u_t^k - \tau \left[ \Kcal_t^* y_{t,1}^{k+1} - \diverg(y_{t,2}^{k+1}) - \diverg(y_{t,3}^{k+1}) - (1-w_t)p_0 \right] + \tau \gamma_t u_{t+1}^k + \tau \gamma_{t-1} u_{t-1}^k}{\tau (\gamma_t + \gamma_{t+1}) +1}$
          \State $z_t^{k+1} - \tau \left[ 2(1-w_t) p_0 + \diverg(y_{t,3}^{k+1}) - \diverg(y_{t,4}^{k+1}) \right]$
          \State {\it Overrelaxation}
          \State $(\bar{u}_t^{k+1}, \bar{z}_t^{k+1}) = 2 (u_t^{k+1}, z_t^{k+1}) - (u_t^k,z_t^k)$
    	\EndFor
	\EndWhile\\
\Return for all $t = 1,\dots,T$: $u_t = u_t^k$
\end{algorithmic}
}
\label{alg:fmri}
\end{algorithm}
To solve the problem, we again perform a proximal gradient descent on the primal variables $\ubold$ and $\zbold$, and a proximal gradient ascent on the dual variables $\ybold$, where the gradients are taken with respect to the linear parts, the proximum with respect to the nonlinear parts. 
We hence need to compute the proximal operators for the nonlinear parts in \eqref{eq:primal_dual} to obtain the update steps for $u_t,z_t$ and $\ybold_t$ for every $t = 1, \dots, T$. 
The proximal operators for $\phi_t(y_{t,1}) = \frac{1}{2 \alpha_t} \| y_{t,1} \|_{\C^{M_t}}^2$ can be computed exactly as in \eqref{eq:prox_dual_l2}.
The proximal operators for the updates of $y_{t,j}$, $j = 2,3,4$, are given by projections onto the sets $C_{t,j}$ similar to \eqref{eq:prox_proj}.
For the squared norm related to the time regularization, we notice that for every $1< t < T$, $u_t$ only interacts with the previous and the following time step, i.e. $u_{t-1}$ and  $u_{t+1}$. 
Hence, analogously to $\phi_t$, the proximum for 
\begin{align*}
	\psi_t(u_t) = \frac{\gamma_{t-1}}{2} \|u_t - u_{t-1}\|_{\C^N}^2 + \frac{\gamma_t}{2} \|u_{t+1} - u_t\|_{\C^N}^2
\end{align*}
is given by 
\begin{align*}
	u_t = \prox_{\tau \psi_t} (r) \Leftrightarrow u_t = \frac{r + \tau \gamma_t u_{t+1} + \tau \gamma_{t-1} u_{t-1}}{\tau (\gamma_t + \gamma_{t-1}) + 1}. 
\end{align*}
The two odd updates for $t = 1$ and $t = T$ can be obtained by the same formula by simply setting $\gamma_0 = 0$ and $\gamma_T = 0$, respectively.
Putting everything together, we obtain Algorithm \ref{alg:fmri}.\\

\subsection{Step sizes and stopping criteria}
We quickly discuss the choice of the step sizes $\tau, \sigma$ and stopping criteria for Algorithm \ref{alg:fmri}. 
In most standard applications it stands to reason to choose the step sizes according to the condition $\tau \sigma \| L \|^2 < 1$ ($L$ denotes the collection of all operators) such that convergence of the algorithm is guaranteed \cite{ChambollePock}.
However, depending on $T$, i.e. the number of time frames we consider, the norm of the operator $L$ 
can be very costly to compute, or too large such that the condition $\tau \sigma \| L \|^2 < 1$ only permits extremely small step sizes. 
For practical use, we instead simply choose $\tau$ and $\sigma$ reasonably ''small`` and track both the energy of the problem and the primal-dual residual \cite{Goldstein:Adaptive} to monitor convergence. 
For the sake of brevity, we do not write down the primal-dual residual for Algorithm \ref{alg:fmri} and instead refer the reader to \cite{Goldstein:Adaptive} for its definition. 
The implementation is then straightforward.
We hence stop the algorithm if both, the relative change in energy between consecutive iterates and the primal-dual residual, have dropped below a certain threshold.

\subsection{Practical considerations}
It is clear that for a large number of time frames $T$ Algorithm \ref{alg:fmri} starts to require an increasing amount of time to return reliable results and for reasonably ''large`` step sizes $\tau$ and $\sigma$ it is even doubtful whether we can obtain convergence. 
In practice, it is hence necessary to divide the time series $\Tbold = \{1, \dots, T\}$ into $l$ smaller bits of consecutive time frames. More precisely, choose numbers $1 \leq T_1 < \dots < T_l = T$ such that $\Tbold = \Tbold_1 \cup \Tbold_2 \cup \dots \cup \Tbold_l$ with $\Tbold = \{1, \dots, T_1, T_1+1, \dots, T_2, \dots, T_{l-1}+1, \dots, T_l\}$.
We can then perform the reconstruction separately for all $\Tbold_i$. 
In order to keep the ''continuity`` between $\Tbold_i$ and $\Tbold_{i+1}$, we can include the last frame of $\Tbold_i$ into the reconstruction of $\Tbold_{i+1}$ by letting $\gamma_{T_i} \neq 0$ and choosing $u_{T_i}$ as the respective last frame of $\Tbold_i$.
This divides the overall problem into smaller and easier subproblems, which can be solved faster.
In practice, we observed that a size of five to ten frames per subset $\Tbold_i$ is a reasonable choice, which essentially gives very similar results as doing a reconstruction for the entire time series $\Tbold$.

\section{Numerical results}\label{sec:numerics}
In this section we show numerical results for simulated fMRI data and real DCE-MRI data. 

\subsection{Methods}
In order to evaluate the potential of the proposed method and to illustrate how each of the terms of our model contributes to the result, we consider the following reconstruction approaches:
\begin{itemize}
\item Least squares (LS) reconstruction of the dynamic sequence: 
\begin{alignat}{4} \label{ls}
\min_{u_t} ~  \| \Kcal_t u_t - f_t \|_{\C^{M_t}}^2   
,\quad t= 1,\ldots,T.
\end{alignat}
The solution of (\ref{ls}) corresponds to a conventional, non-regularized reconstruction.

\item Reconstruction of the dynamic sequence using spatial TV regularization (TV):
\begin{alignat}{4} \label{tv}
\min_{u_t} ~ \frac{\alpha_t}{2} \| \Kcal_t u_t - f_t \|_{\C^{M_t}}^2 + \TV(u_t)  
,\quad t= 1,\ldots,T
\end{alignat}
The solution of (\ref{tv}) gives a gradient sparsity promoting reconstruction without time correlation between the frames $u_t$ in the image sequence.

\item Reconstruction of the dynamic sequence using temporal (smoothness) regularization, but no spatial regularization (temp):
\begin{alignat}{4} \label{l2t}
\min_{\ubold} ~ \sum_{t=1}^T \frac{\alpha_t}{2} \| \Kcal_t u_t - f_t \|_{\C^{M_t}}^2  + \sum_{t=1}^{T-1} \frac{\gamma_t}{2} \|u_{t+1} - u_t \|_{\C_N}^2 
\end{alignat}

\item Reconstruction of the dynamic sequence using spatial TV regularization and temporal regularization (temp + TV):
\begin{alignat}{4} \label{tvt}
\min_{\ubold} ~ \sum_{t=1}^T \frac{\alpha_t}{2} \| \Kcal_t u_t - f_t \|_{\C^{M_t}}^2 + \sum_{t=1}^T \TV(u_t)  + \sum_{t=1}^{T-1} \frac{\gamma_t}{2} \|u_{t+1} - u_t \|_{\C_N}^2 
\end{alignat}
The solution of (\ref{tvt}) is similar to the spatio-temporal regularization techniques in \cite{adluru2009,wundrak2016} and provides a reference
of spatio-temporal regularization without the structural prior.

\item The proposed approach using temporal regularization and $\ICBTV$ regularization for the utilization of the structural prior information:
\begin{alignat}{4} \label{icbtv}
\min_{\ubold} ~  &\sum_{t=1}^T \frac{\alpha_t}{2} \| \Kcal_t u_t - f_t \|_{\C^{M_t}}^2 \notag\\
+ &\sum_{t=1}^T \left(w_t \TV(u_t) + (1-w_t) \ICBTV^{p_0}(u_t,u_0) \right)  \notag\\
+ &\sum_{t=1}^{T-1} \frac{\gamma_t}{2} \|u_{t+1} - u_t \|_{\C_N}^2 
\end{alignat}
\end{itemize}
In the following we show results for two test cases on simulated and experimental dynamic MRI data, respectively.

\subsection{Simulated data}
\label{subsec:artifical data}
In this paragraph, we first describe the generation of a synthetic test data set from a human brain phantom in the style of fMRI measurements. Afterwards we give various illustrations of the corresponding results obtained with the aforementioned reconstruction methods.  

\subsubsection{Generation of a simulated data set}
As described before, in fMRI data is typically acquired in two stages: after a high-resolution scan that serves as anatomical reference the actual dynamic sequence is sampled. 
While conventionally the former is measured in T$1$ contrast, the changes in the cerebral blod flow are detectable in a T$2^*$ weighted MR signal such that the dynamic sequence needs to be acquired with this contrast. 
In particular, one has to face different image contrasts in the anatomical prior and the dynamic sequence.

To mimic such fMRI measurements, we used a slice of a realistically simulated MR phantom from the BrainWeb database \cite{cocosco:brainweb}, consisting of T$1$ and corresponding T$2$ axial MR images of the human brain of a size of $109 \times 91$ pixels each, with pixel intensities in $[0,1]$. 
To create a dynamic ground truth sequence, we duplicated the T$2$ image to obtain a sequence of 60 time steps $t$. 
We then chose a small area of gray matter (see the red box in Figure \ref{fig:results_synth}) and changed the intensity in this area according to a canonical hemodynamic response function (HRF) from the SPM software package\footnote{\texttt{http://www.fil.ion.ucl.ac.uk/spm/}}.
The HRF (see ``ground truth" in Figure \ref{fig:timeplots_synth} and \ref{fig:timeplots_single}) consists of a linear combination of two Gamma functions, where we used the default parameters of the package.
The amplitude of the HRS function has been scaled to approximately $0.1$, which corresponds to $10$ per cent of the overall image scale. 
It is worth noticing that this change is very small in comparison to some of the image intensities. 
This can e.g. be observed in the last line of Figure \ref{fig:results_zoom}, which shows a zoom into the red box in Figure \ref{fig:results_synth}.

To obtain a data set, we applied the forward operator \eqref{eq:forward_op_art} to this synthetic ground truth.
For the anatomical T$1$ prior, we set $\Pcal_0$ to be the identity leading to a full Cartesian sampling of $109 \times 91$ Fourier coefficients in $k$-space. 
For the dynamic T$2$ sequence, we chose $\Pcal_t$ such that at each time step frequencies located on 5 discrete spokes through the center of the $k$-space are measured (see first line in Figure \ref{fig:results_synth}). 
The spokes angles are chosen according to golden angle sampling (cf. e.g. \cite{winkelmann2007}).
This sampling scheme amounts to measuring less than $5 \%$ of the $109 \times 91$ Fourier coefficients of the full Cartesian sampling at each time step.
We finally added noise of approximately $5$ per cent of the energy of the whole signal to the sampled Fourier frequencies.
 
\subsubsection{Reconstructions}
\begin{figure}[ht!]
	\includegraphics[width=0.95\textwidth]{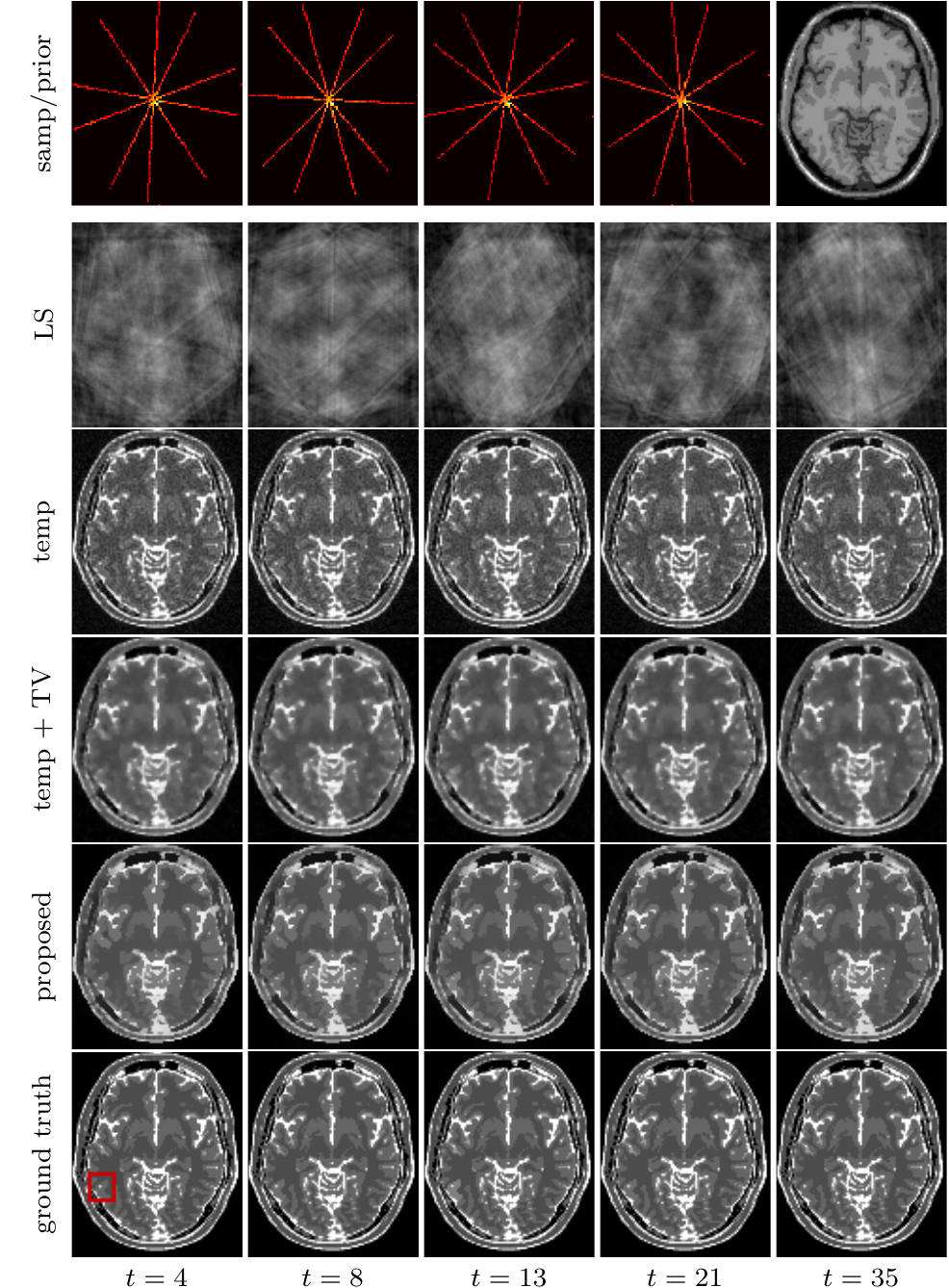}
    \caption{Reconstructions using the different methods (LS), (temp), (temp + TV) and the proposed method at five different time points.
    The first line shows the respective $k$-space samplings and the anatomical prior used for the proposed method.
    The red box marks the activated region of interest, where the signal changes are observable.
    Figure \ref{fig:results_zoom} provides a zoom into this region.}
    \label{fig:results_synth}
\end{figure}
\begin{figure}[ht!]
	\includegraphics[width=\textwidth]{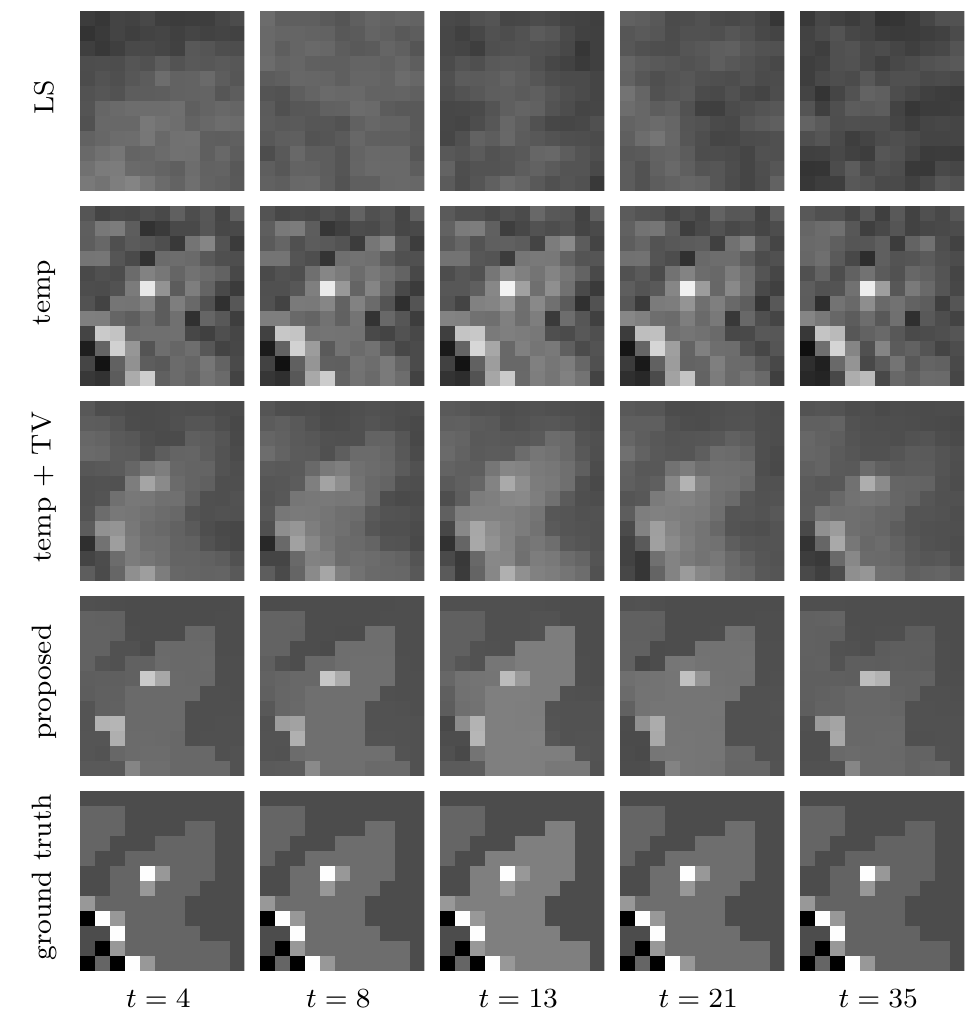}
    \caption{Zoom into the region of interest marked with the red box in Figure \ref{fig:results_synth} for the different methods (LS), (temp), (temp + TV) and the proposed method.
    The ground truth on the bottom line shows the true activation of the gray matter.}
    \label{fig:results_zoom}
\end{figure}
\begin{figure}[ht!]
	\includegraphics[width=\textwidth]{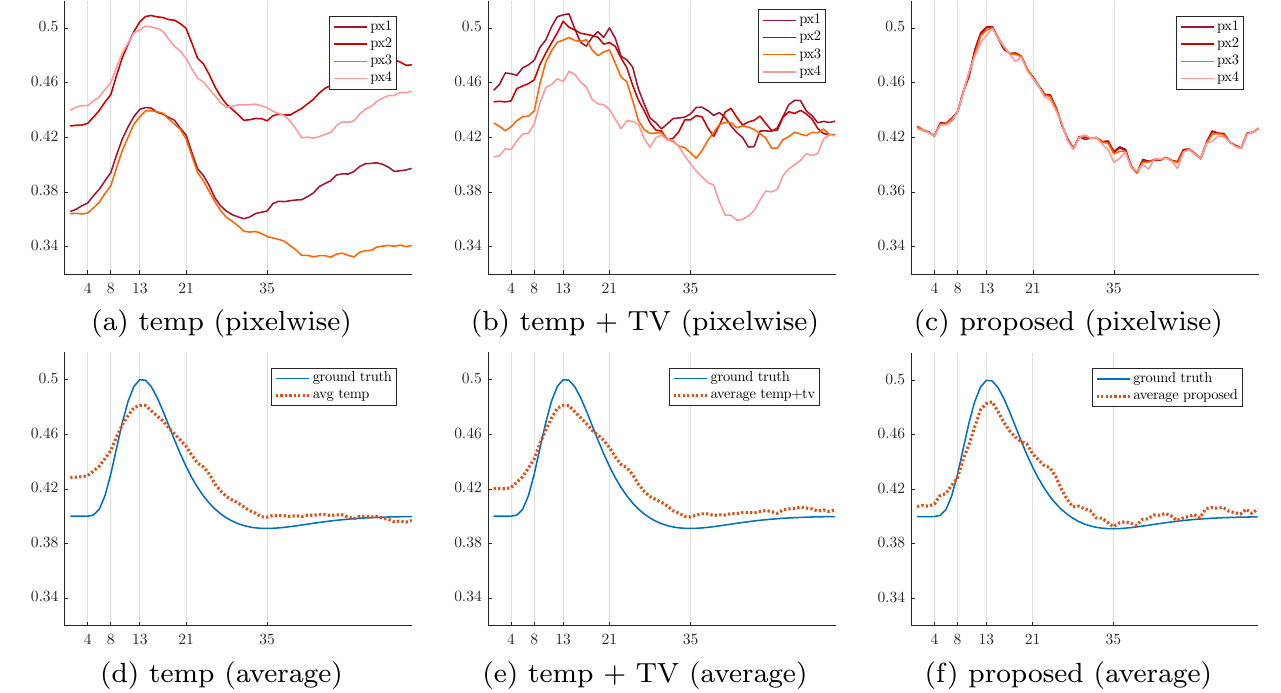}
    \caption{Ground truth and reconstructions of the (simulated) hemodynamic response over 60 time frames for the different methods (temp), (temp + TV) and the proposed method: (a)-(c) signal from four adjacent single pixels from the region of interest over time, (d)-(f) average over the entire activated region over time. 
    }
    \label{fig:timeplots_single}
\end{figure}
\begin{figure}[t]
	\includegraphics[width=\textwidth]{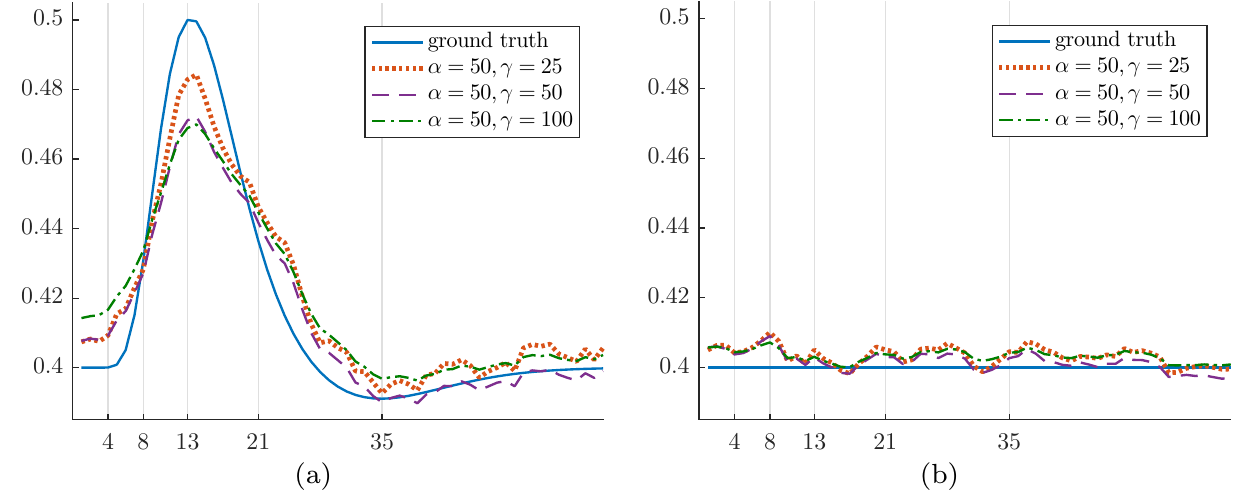}
    \caption{Ground truth and reconstructions with the proposed method of the (simulated) hemodynamic response over 60 time frames for different amounts of temporal regularization.
    (a) Average of the reconstructed hemodynamic response over the whole activated region.
    (b) Average over an inactive region.}
    \label{fig:timeplots_synth}
\end{figure}

Figure \ref{fig:results_synth} shows results for the different reconstruction methods for five different times $t= 4,8,13,21,35$, together with the $k$-space sampling and the anatomical prior. 
The prior has been reconstructed using \eqref{tv} and has then been used to create an artificial subgradient $q_0^\eta$ as in \eqref{eq:subgrad_eta} with threshold $\eta = 0.05$.
Note that we did not include single-frame TV reconstructions \eqref{tv}, since the image quality did not increase significantly over the quality of a least squares reconstruction \eqref{ls}.
The presented times $t$ have been chosen to show different states of activation of the region of interest (ROI), which is marked by a red box in the bottom line of Figure \ref{fig:results_synth}. 
Figure \ref{fig:results_zoom} provides a close-up into the ROI.
For $t=4$ there is no activation, $t=8$ shows half of the increase of the signal, $t=13$ shows maximum activation, $t=21$ shows half of the decrease of the signal and $t=35$ a slight undershoot after the signal has decayed.
These times are also marked in Figures \ref{fig:timeplots_synth} and \ref{fig:timeplots_single}. 
The former shows the average signal over time over the activated region in the ROI in comparison to a non-activated control region for different choices of regularization parameters $\alpha$ and $\gamma$.
The latter figure shows the average signal over time over the activated region in the ROI for the different methods \eqref{ls}, \eqref{l2t}, \eqref{tvt} and \eqref{icbtv}, together with the signal from four individual adjacent pixels from the activated region.   
The parameters used to obtain the results are $\gamma_t = \alpha_t = 1$ for (temp) , $\alpha_t = \gamma_t = 500$ for (temp + TV) and $\alpha_t = 50$, $\gamma_t = 25$ for the proposed method. 
Note that we keep the parameters constant for all $t = 1, \dots, T$.

In Figure \ref{fig:results_synth} we observe that while a simple single-frame least squares reconstruction \eqref{ls} does not yield satisfactory results, a temporal regularization \eqref{l2t}, resulting in a pixelwise regularity over time, can improve image quality by far. 
However, while the structures of the brain are clearly visible and also the activation in the ROI, the quality of the reconstructions is degraded by heavy noise. 
This is also confirmed by the plots in Figure \ref{fig:timeplots_single}(a) and \ref{fig:timeplots_single}(d). 
While the reconstruction of the average signal over the activated region is decent (\ref{fig:timeplots_single}(d)), single pixels from the activated region show a lot of variance (\ref{fig:timeplots_single}(d)). 
The main reason for this is the noise in the data and the lack of any spatial regularity.

Adding a spatial TV regularization to the temporal regularization, we obtain method \eqref{tvt}, which leads to an overall reduction of noise and more spatial regularity in the reconstructed time frames, but also blurs out some of the structures.  
The activated region is as well visible, and the averaged signal over the ROI gives a similar result than before, which can be observed in Figure \ref{fig:timeplots_single}(e). 
However, we as well observe a greater coherence among the adjacent single pixels from the ROI depicted in Figure \ref{fig:timeplots_single}(b). 

The reconstructions using the proposed method including the anatomical prior show the most convincing results. 
In particular, the use of the anatomical prior leads to a very high spatial regularity with distinct edges such that most of the features of the brain are clearly visible. 
Furthermore, both the average over the activated region {\it and} the single pixels now show almost the same signal, which comes closest to our simulated ground truth (see Figure \ref{fig:timeplots_single}(c) and (f)). 
We like to point out that the reconstruction also shows the slight undershoot after the decay of the signal at $t=35$, which can be seen in Figure \ref{fig:timeplots_synth}.

We also comment on a few imperfections. 
All of the reconstructed signals show a slight bias in the sense that the amplitude of the signal has been slightly decreased. 
This is due to the (necessary) temporal regularization. 
The problem is that the amplitude of the signal varies over time, which would actually require a lower time regularity where the change between consecutive time frames is high, and a higher regularity where the difference is small. 
This can be observed in Figure \ref{fig:timeplots_synth}, where a higher $\gamma$ decreases the amplitude of the signal while enforcing less oscillations of the signal in the flat region after $t=35$.
Hence the choice of $\gamma$ always results in a trade-off between temporal smoothness and amplitude of the signal. 
If one has a rough estimate on when the signal develops, the proposed method offers a remedy by choosing a smaller $\gamma_t$ while it develops, and a larger $\gamma_t$ in the flat region at the end instead of a global $\gamma$.
This could further increase the quality of the results. 

\subsection{Experimental data from small animal imaging}
\label{experes}
Subsequently, we first describe the preparation and execution of an experiment to generate real small animal DCE-MRI data and give details on the acquisition protocol.
Afterwards, we show various visualizations of the respective results obtained with the aforementioned reconstruction methods.

\subsubsection{Animal preparation}
All animal experiments were approved by the Animal Health Welfare and Ethics Committee of University of Eastern Finland. 
1x106 C6 (ECACC 92090409) rat glioma cells (Sigma) were implanted into the brain of a 200g female Wistar rat under ketamin/medetomidine hydrochloride anesthesia. 
Tumor imaging was performed 10 days post-implantation. 
During the experiments, the animal was anesthetized with isoflurane (5\% induction, 1-2\% upkeep) and kept in fixed position in a holder which is inserted into the magnet. 
A needle was placed into the tail vein of the animal for the injection of the contrast agent.

\subsubsection{Acquisition of the data}
All MR data were collected using a 9.4 T horizontal magnet interfaced to Agilent imaging console and a volume coil transmit/quadrature surface coil receive pair (Rapid Biomed, Rimpar, Germany). 
For the dynamical data a gradient-echo based radial pulse sequence 
with repetition time 38.5 ms, echo time 9 ms, flip angle 30 degrees, field-of view 32x32 mm, slice thickness 1.5 mm, number of points in each spoke 128, 610 spokes collected in sequential order with an golden angle interval of 111.246 degrees before repeating the same cycle of spokes 
for 25 times, leading to overall measurement sequence of 15250 spokes of data. 
Measurement time for a full cycle of 610 spokes was 
$610 \cdot 38.5 {\rm ms} = 23.46 {\rm s}$. 
Gadovist (1mmol/kg) was injected i.v. after one minute from the beginning of the dynamic scan over a period of 3s. 

Anatomical reference images were acquired from the same slice before and after the dynamical experiment using a gradient-echo pulse sequence with repetition time of 1s and echo time of 2.8ms
with Cartesian sampling of 128x128 points of k-space data.

\subsubsection{Reconstructions}

\begin{figure}[ht!]
\centerline{\includegraphics[width=0.8\textwidth]{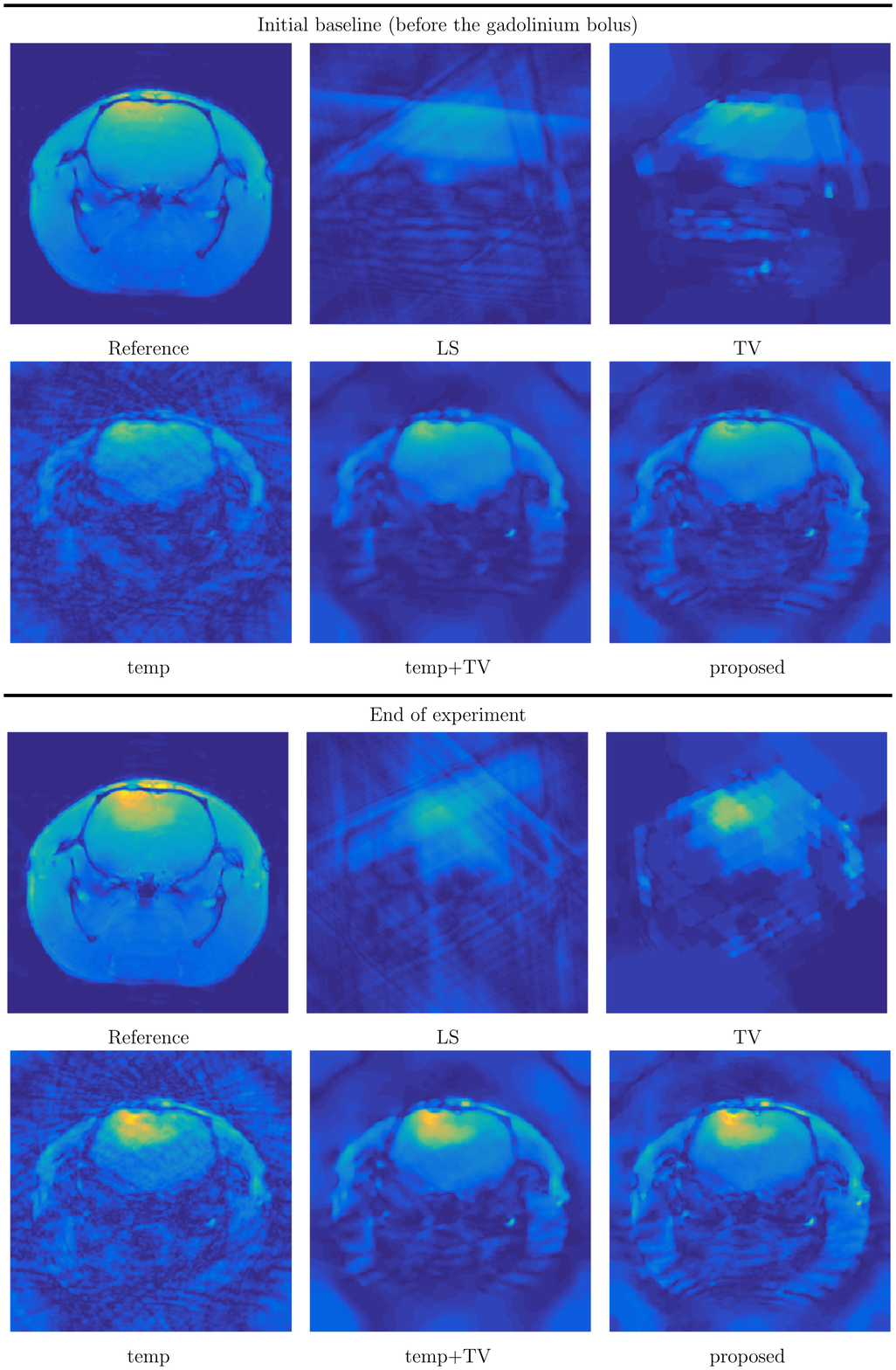}}
\caption{DCE-MRI using radial golden angle data from a glioma model rat specimen. Top half: Anatomical reference image before the dynamical experiment and reconstructed frame
from the initial base line measurement at $t=200$. 
Top row shows the reference, the LS estimate and the estimate with the TV regularization. The second row shows the reconstructions using temporal regularization, TV combined with temporal regularization and the $\ICBTV$ combined with temporal regularization.
The bottom half of image shows in respective order the anatomical reference after the experiment and the last frame of the reconstructed time series.
The number of radial spokes for each time frame was three.}
\label{efig1}
\end{figure}

\begin{figure}[ht!]
 \centerline{\includegraphics[width=0.85\textwidth]{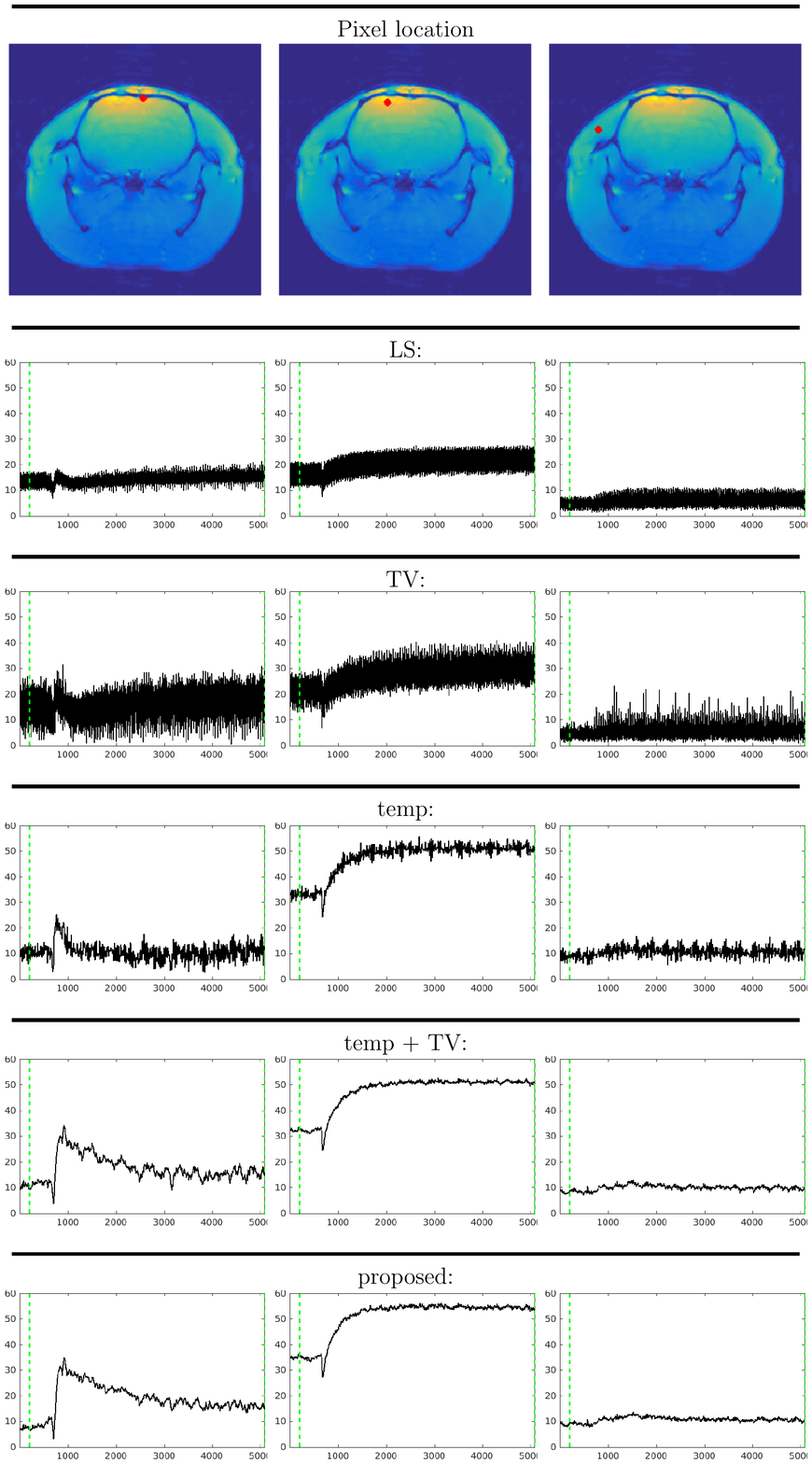}}
\caption{DCE-MRI using radial golden angle data from a glioma model rat specimen. Top: Anatomical reference image before the dynamical experiment. The squares mark the
location of the pixel of interest for each column. 
Rows two to six show time series of the image magnitude in the selected pixels using the LS, TV regularization, temporal  regularization, TV regularization combined with temporal  regularization and $\ICBTV$ combined with temporal regularization, respectively. 
The dashed vertical lines denote the time points of the slices shown in Figure \ref{efig1}.}
\label{efig2}
\end{figure}

\begin{figure}[ht!]
\centerline{\includegraphics[width=0.82\textwidth]{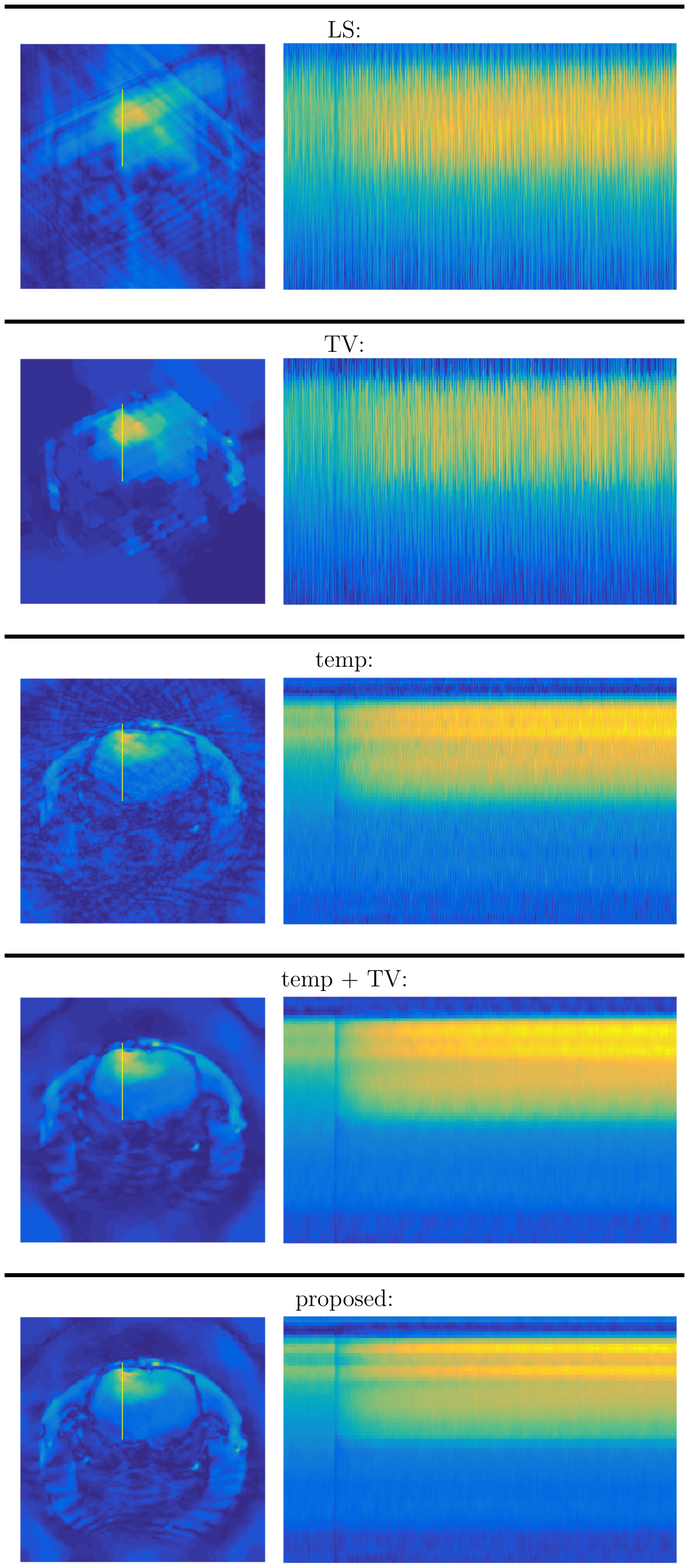}}
\caption{DCE-MRI using radial golden angle data from a glioma model rat specimen. The time series of pixel values in a vertical ROI line through the brain and the glioma. The 
location of the ROI line is indicated with a yellow line in the reconstructions of the last frame on the left. The time series of the pixel values are shown in the right column such that horizontal axis corresponds to time and vertical to location along the ROI line. Rows from top to bottom: 
LS, TV regularization, temporal regularization, TV combined with temporal regularization and $\ICBTV$ combined with  temporal regularization, respectively.}
\label{efig3}
\end{figure}

The prior image $u_0 \in \C^{65536}$ for the structural regularization was 
obtained as solution of (\ref{tvu0}) using the Cartesian reference data measured before the dynamical experiment. 
The value of the weighting parameter in (\ref{tvu0}) was $\alpha_0 = 2$. 
Note that due to the Cartesian sampling in $k$-space, the prior could also be obtained using an inverse Fourier transform, however yielding a $128 \times 128$ resolution in this case. 
Using method \eqref{tvu0} and a finer grid of $256 \times 256$ enables us to ``lift'' the resolution of the prior, which can then also be transferred to the dynamic sequence.
The magnitude of the prior image $u_0$ is shown in the first image of the top row in Figure \ref{efig1}. 

For the dynamical reconstructions, the dynamical data was divided into 
$T=5083$ time steps, 
each set consisting of three consecutive spokes, leading to measurement data $f_t \in \C^{384}$ for each frame in the time series.  
As for the prior, the unknown images were 
represented using a $256 \times 256$ pixel grid for the reconstructions, leading to unknown $u_t \in \C^{65536}$ for each time instant $t$. 
The parameters for the reconstructions were selected manually by pilot runs. The value  
$\alpha_t = 5$ was used in all the regularized solutions (\ref{tv}-\ref{icbtv}) and the value of the temporal regularization parameter was $\gamma_t = 20$ in the models (\ref{l2t}-\ref{icbtv}) containing temporal regularization. The value of the weighting parameter for the structural regularization in (\ref{icbtv}) was selected as $w_t = 0.1$. 

The results using the experimental data are shown in Figures \ref{efig1}-\ref{efig3}. The top half of Figure \ref{efig1} shows the anatomical reference before the experiment and magnitude of reconstructed frames $u_t$ at $t=200$ from the initial baseline measurement before the injection of the gadolinium bolus. 
The first row shows the anatomical reference, the LS reconstruction (\ref{ls}) and the reconstruction with the TV regularization (\ref{tv}). 
The second row shows the reconstruction (\ref{l2t}) using temporal regularization, (\ref{tvt}) with combination of temporal regularization and TV, and the proposed approach (\ref{icbtv}). 
The bottom half of Figure \ref{efig1} shows in the same order the anatomical reference after the experiment and the last frames of the reconstructed time series.    
Notice that the muscle tissues are present more clearly in the reference images than in the dynamical images due to the different echo and repetition times in the dynamical and reference measurements. 

Figure \ref{efig2} shows the time-series of the magnitude of reconstructed pixel value at three selected pixels. Each column shows the data for a single pixel, which is indicated by a red dot mark in the anatomical reference shown in the top row. 
The pixel time series in rows two to five are LS (\ref{ls}), TV regularization (\ref{tv}), temporal regularization (\ref{l2t}), temporal regularization and TV (\ref{tvt}), and the proposed method (\ref{icbtv}).  
The dashed vertical lines in the time series plots mark the time locations of the reconstructed frames in Figure \ref{efig1}.
Figure \ref{efig3} shows time series of reconstructed pixel values in a vertical ROI 
line through the brain and the glioma.
The estimates are from top to bottom in the same order as in Figure \ref{efig2}.
The vertical ROI line is indicated in the reconstructions
of the last time frame on the left.
The time series of the reconstructed ROI values are shown as intensity images on the right such that the horizontal axis corresponds to frame number and vertical axis to pixel location along the vertical ROI line.

The findings from the results with the experimental DCE data are very similar to those with the simulated fMRI data. First, we observe that the single frame reconstructions, namely the LS and TV, do not yield satisfactory results; the images are spatially of poor quality and the recovered time series are very noisy, rendering reliable visual or quantitative analysis of DCE parameters infeasible.
Using the temporal regularization (\ref{l2t}), improves the results clearly over the uncorrelated reconstructions. 
The images show well the structure of the brain and also the dynamics caused by the gadolinium contrast agent are clearly visually detectable, but the results are still contaminated with a high level of noise. 

Using the combination of temporal and TV regularization (\ref{tvt}), improves the results significantly, leading to a clear reduction of noise in the reconstructed images and time series. 
However, similarly as in the simulated test case, this comes with the cost of some blurring of the images, especially on the tissue interfaces.
The proposed method (\ref{icbtv}), which combines the time regularization and structural prior, produces again the most convincing results. 
While the reconstructed pixel time series are visually not much different from the results with (\ref{tvt}), the spatial features of the reconstructions are sharper.

\section{Conclusion and outlook}
\label{sec:outlook}
In the previous sections we presented a novel variational model for the reconstruction of highly subsampled dynamic MRI data where an anatomical scan (at high spatial resolution) has been acquired prior to the dynamic sequence. 
Combining radial golden angle sampling with a suitable time regularization, spatial TV regularization and with the infimal convolution of TV Bregman distances allowing to incorporate the structural information of the anatomical prior, we obtained spatially highly resolved reconstructions at a high temporal resolution. 

Summing up the results of tests on a simulated data set based on fMRI as well as on experimental small animal DCE-MRI data, we draw the following conclusions:
naturally, a simple least squares (LS) reconstruction of each individual frame could not provide meaningful results due to the severe undersampling. 
Adding a spatial TV regularization did not significantly increase the quality of the reconstructed images. 
As expected, this approach yielded piecewise constant images, but the ratio of sampled Fourier coefficients in $k$-space in comparison to the desired spatial resolution of the reconstructions was too small to obtain reasonable results. 
Remarkably, adding only Tikhonov regularization on the time derivative without any additional spatial regularization already resulted in by far more meaningful reconstructions than the LS approach, while the obtained images were still corrupted by heavy noise.
Integrating spatial TV regularization to the aforementioned model removed most of the noise and indeed provided high quality reconstructions. 
Incorporating structural information from the anatomical prior, we could then obtain very detailed results despite the severe subsampling enabling high temporal resolution.

In view of these promising results, we state some open questions and sketch additional ideas whose detailed study is left to future research.

We used the infimal convolution of TV Bregman distances to incorporate structural information from the anatomical prescan. 
Naturally, this gives rise to the question whether alternative means of incorporating structural prior information such as the concepts of weighted total variation (wTV) or directional total variation (dTV), respectively, both proposed in \cite{Ehrhardt2016}, yield significantly different results. 
In any case, it would be interesting to see how such a modified approach compares to the method proposed in this paper concerning quality of the reconstructed images, but also regarding computational complexity of solving the respective minimization problem. 

Moreover, the temporal coupling of time frames serves as a further starting point for future research. 
Here, we decided to apply Tikhonov regularization of the time derivative, however, one could also argue in favor of other concepts: 
since in the areas of application we considered in this paper the dynamic changes happen to take place in only a small portion of the entire image domain, decomposition of the dynamic sequence into a low rank part $L$ and a part $S$ which is sparse in some transform domain \cite{Tremoulheac:lowRankPlusSparsePrior,Otazo:lowRankPlusSparseMatrixDecomposition} could be an interesting alternative. 
Assuming that the dynamic changes mainly are contained in $S$, while $L$ ideally comprises the part staying constant over time, it seems particularly reasonable that the structures of the constant part of every time frame bear close resemblance to the structure of the anatomical prior image. 
Hence it would stand to reason to apply the infimal convolution of TV Bregman distances only to the low rank part leaving the sparse part untouched. 
However, against the backdrop of different dimensions of the low rank part of the dynamic sequence $L$ and the anatomical prior image $u_0$ it is not yet clear what would be the most suitable way of solving the corresponding optimization problem. 

Finally, with respect to experimental data, a more careful correction of artifacts due to different acquisition protocols between anatomical prescan and the dynamic sequence might be an interesting aspect.

\section*{Acknowledgements}

This work has been supported  ERC via Grant EU FP 7 - ERC Consolidator Grant 615216 LifeInverse, by the German Ministery for Science and Education (BMBF) through the project MED4D, by Academy of Finland (Finnish Programme for Center of Excellence in Research 2012-2017, project 250215) and Jane and Aatos Erkko Foundation. 
The authors would like to thank the Isaac Newton Institute for   Mathematical Sciences, Cambridge, for support and hospitality during  the programme Variational Methods for Imaging and Vision, where work on this paper was  undertaken, supported by EPSRC grant no EP/K032208/1.  

\bibliographystyle{abbrv}
\bibliography{fmri_ref}

\end{document}